\documentclass[12pt]{article}
\usepackage[a4paper,left=2cm,right=1.2cm,top=3cm,bottom=4cm]
{geometry}

\usepackage{amsmath,amstext, amsthm, empheq, extpfeil,  
amsbsy,amssymb,marvosym,fancyhdr,graphicx,amscd,amsfonts,latexsym,delarray,stackrel,lineno,color,cite,appendix,xcolor,url,hyperref}
\usepackage{wasysym}
\usepackage{mdframed}

\usepackage{srcltx}
\usepackage{mathrsfs}
\usepackage{float} 
\usepackage{subfig} 
\usepackage[all]{xy}
\usepackage{t1enc}
\usepackage{mathrsfs}
\usepackage{pifont}

\usepackage{mathpple}
\usepackage[T1]{fontenc}

\newcommand*\patchAmsMathEnvironmentForLineno[1]{%
  \expandafter\let\csname old#1\expandafter\endcsname\csname #1\endcsname
  \expandafter\let\csname oldend#1\expandafter\endcsname\csname end#1\endcsname
  \renewenvironment{#1}%
     {\linenomath\csname old#1\endcsname}%
     {\csname oldend#1\endcsname\endlinenomath}}%
\newcommand*\patchBothAmsMathEnvironmentsForLineno[1]{%
  \patchAmsMathEnvironmentForLineno{#1}%
  \patchAmsMathEnvironmentForLineno{#1*}}%
\AtBeginDocument{%
\patchBothAmsMathEnvironmentsForLineno{equation}%
\patchBothAmsMathEnvironmentsForLineno{align}%
\patchBothAmsMathEnvironmentsForLineno{flalign}%
\patchBothAmsMathEnvironmentsForLineno{alignat}%
\patchBothAmsMathEnvironmentsForLineno{gather}%
\patchBothAmsMathEnvironmentsForLineno{multline}%
}

\definecolor{Green}{rgb}{0,1,0}
\definecolor{Blue}{RGB}{0,0,191}
\definecolor{mathmodecolor}{RGB}{0,102,0}
\definecolor{keywordcolor}{RGB}{0,51,151}
\definecolor{sourcebackgroundcolor}{RGB}{255,247,223}
\definecolor{unixagred}{RGB}{255,0,0}
\definecolor{lightgray}{RGB}{191,191,191}
\definecolor{green}{RGB}{1,191,191}

\def\wnt{{\widehat{\N^{\times}}}}

\def\cutint{{\int \!\!\!\!\!\! -}}

\def\F{{\mathbb F}}

\def\qqq{\,,\quad \forall}

\def\pert{{\rm Pert}}

\def\bm2{{\rm \B mod^2}}
\def\b2{{\rm \B mod^{\mathfrak s}}}
\def\hatz{{\hat\Z^*}}

\def\Aut{{\rm Aut}}

\def\Ext{{\rm Ext}}

\def\Ker{{\rm Ker}}

\def\Tr{{\rm Tr}}

\def\A{{\mathbb A}}
\def\C{{\mathbb C}}
\def\F{{\mathbb F}}

\def\N{{\mathbb N}}

\def\Q{{\mathbb Q}}
\def\R{{\mathbb R}}
\def\T{{\mathbb T}}
\def\Z{{\mathbb Z}}
\def\H{{\mathbb H}}

\def\B{{\mathbb B}}

\def\Tr{{\rm Tr}}

\def\cA{{\mathcal A}}

\def\cC{{\mathcal C}}

\def\cH{{\mathcal H}}

\def\cK{{\mathcal K}}

\def\cM{{\mathcal M}}

\def\cS{{\mathcal S}}

\def\vp{\varphi}

\def\part{\partial}

\newcommand{\ie}{{\it i.e.\/}\ }
\newcommand{\eg}{{\it e.g.\/}\ }
\newcommand{\cf}{{\it cf.\/}\ }
\newcommand{\opcit}{{\it op.cit.\/}\ }

\def\H{{\mathbb H}}

\def\elel{(  \downarrow 1)}
\def\dodo{(   \downarrow 3)}

\def\rmax{\R_+^{\rm max}}

\def\fr{{\rm Fr}}

\def\rmax{\R_+^{\rm max}}

\def\higgs{{\bf H}}
\definecolor{trust}{rgb}{0,1,1}

\def\nt{\N^{\times}}

\def\lapa{\triangle}

\def\modu{\Delta}

\def\perta{\triangle_\vp}

\def\pertc{\triangle^{(0, 1)}_\vp}

\def\Ker{{\rm Ker}}

\def\noncommutative geometry{{noncommutative geometry }}

\def\qqq{\,,\quad \forall}

\title
{Noncommutative Geometry, the spectral standpoint}

\author{Alain Connes}

\begin{document}

\maketitle

\centerline{\em{In memory of John Roe, and in recognition of his pioneering achievements }}
\centerline{\em{ in coarse geometry and index theory.}}

\begin{abstract} We update our Year 2000 account of Noncommutative Geometry in \cite{ncg2000}. There, we described the following basic features of the subject:

$\blacktriangleright$~ The natural ``time evolution" that makes noncommutative spaces dynamical from their measure theory.

$\blacktriangleright$~ The new calculus which is based on operators in Hilbert space, the Dixmier trace and the Wodzicki residue.
 
$\blacktriangleright$~ The spectral geometric paradigm which 	extends the Riemannian paradigm to the noncommutative world and gives a new tool to understand the forces of nature as pure gravity on a more involved geometric structure mixture of continuum with the discrete.
 
$\blacktriangleright$~ The key examples such as duals of discrete groups, leaf spaces of foliations and deformations of ordinary spaces, which showed, very early on, that the relevance of this new paradigm was going far beyond the framework of Riemannian geometry.
 
 Here, we shall report\footnote{this report makes no claim to be exhaustive and we refer the reader to the collection of surveys \cite{chamssuijl,CCsurvey,DS5,Dundas, GE,FK3,GJV,LM1,LSZ1,Vo1,Wsurvey,XY2},  published in {\em  Advances in Noncommutative Geometry}, Springer,  ISBN 978-3-030-29596-7,
(2019), for a more complete account and bibliography.} on the following highlights from among the many discoveries made since 2000: 
 \begin{enumerate}
 \item The interplay of the geometry with the modular theory for noncommutative tori.
 \item Great advances  on the Baum-Connes conjecture, on coarse geometry and on higher index theory.
 \item The geometrization of the pseudo-differential calculi using smooth groupoids.
 \item The development of Hopf cyclic cohomology. 
 \item The increasing role of topological cyclic homology in number theory, and of the lambda operations in archimedean cohomology. 
 \item The understanding of the renormalization group as a motivic Galois group.
 \item The development of quantum  field theory on noncommutative spaces.
 \item The discovery 
  of a simple equation whose irreducible representations correspond to $4$ dimensional spin geometries with quantized volume and give an explanation of the Lagrangian of the standard model coupled to gravity.
 \item The discovery that very natural toposes such as the scaling site provide the missing algebro-geometric structure on the noncommutative adele class space underlying the spectral realization of zeros of $L$-functions.

\end{enumerate}
\end{abstract}


\tableofcontents

\section{Introduction, two key examples}

Noncommutative geometry has its roots both in quantum physics and in pure mathematics. One of its main themes is to explore and investigate ``new spaces" which lie beyond the scope of standard tools. In physics the discovery of such new spaces goes back to Heisenberg's matrix mechanics which revealed the non-commutative nature of the phase space of quantum systems. The appearance of a whole class of ``new spaces" in mathematics  and the need to treat them with new tools came later. It is worth explaining in simple terms when one can recognize that a given space such as the space of leaves of a foliation, the space of irreducible representations of a discrete group, the space of Penrose tilings, the space of lines\footnote{infinite geodesics} in a negatively curved compact Riemann surface, the space of rank one subgroups of $\R$, etc, etc, is ``noncommutative" in the sense that classical tools are inoperative. The simple characteristic  of such spaces is visible  even at the level of the underlying ``sets": for such spaces their cardinality is the same as the continuum but nevertheless it is not possible to put them constructively in bijection with the continuum. More precisely any explicitly constructed map from such a set to the real line fails to be injective! It should then be pretty obvious to the reader that there is a real problem when one wants to treat such spaces in the usual manner. One may of course dismiss them as ``pathological" but this is ignoring their abundance since they appear typically when one defines a space as  an inductive limit. Moreover extremely simple toposes admit such spaces as their spaces of points and discarding them is an act of blindness. The reason to call such spaces ``noncommutative" is that if one accepts to use noncommuting coordinates to encode them, and one extends the traditional tools to this larger framework, everything falls in place. The basic principle is that   one should take advantage of the  presentation of the space as a quotient of an ordinary space by an equivalence relation and instead of effecting the quotient in one stroke, one should associate to the equivalence relation its convolution algebra over the complex numbers, thus keeping track of the equivalence relation itself. The noncommutativity of the algebra betrays the identification of points\footnote{The groupoid identifying two points gives the noncommutative algebra of matrices $M_2(\C)$} of the ordinary space and when the quotient happens to exist as an ordinary space, the algebra is Morita equivalent  to the commutative algebra of complex valued functions on the quotient. One needs of course to refine the perception of the space one is handling according to the natural hierarchy of geometric points of view, going from measure theory, topology, differential geometry to metric geometry.

It is important to stress from the start that noncommutative geometry is not just a ``generalization" of ordinary geometry, inasmuch as while it covers as a special case the known spaces, it handles them from a completely different perspective, which is inspired by quantum mechanics rather than by its classical limit. Moreover there are many surprises and the new spaces admit  features which have no counterpart in the ordinary case, such as their time evolution which appears already from their measure theory and turns them into thermodynamical objects.

 The obvious question  then, from a conservative standpoint,  is why one needs to explore such new mathematical entities, rather than staying in the well paved classical realm ? The simple answer is that noncommutative geometry allows one to comprehend in a completely new way the following two spaces whose  relevance can hardly be questioned :
 
 \newpage


\begin{mdframed}[backgroundcolor=yellow!2]
\centerline{\underline{\bf Space-Time  and NCG}}

\vspace{0.5cm}

The simplest reason why noncommutative geometry is relevant for the understanding of the geometry of space-time is the key role of the non-abelian gauge theories in the Standard Model of particles and forces. The gauge theories enhance the symmetry group of gravity, \ie the group of diffeomorphisms ${\rm Diff}(M)$ of space-time, to a larger group which is a semi-direct product with the group of gauge transformations of second kind. Searching for a geometric interpretation of the larger group as a group of diffeomorphisms of an higher dimensional manifold is the essence of the Kaluza-Klein idea. Noncommutative geometry gives another track. Indeed the group of diffeomorphisms ${\rm Diff}(M)$ is the group of automorphisms of the star algebra $\cA$ of smooth functions on $M$, and if one replaces $\cA$ by the (Morita equivalent) noncommutative algebra $M_n(\cA)$ of matrices over $\cA$ one enhances the group of automorphisms of the algebra in exactly the way required by the non-abelian gauge theory with gauge group ${\rm SU}(n)$. The Riemannian geometric paradigm is extended to the noncommutative world in an operator theoretic and spectral manner. A geometric space is encoded by its algebra of coordinates $\cA$ and its ``line element" which specifies the metric. The new geometric paradigm of spectral triples (see \S \ref{sectspecgeom}) encodes the discrete and the continuum on the same stage which is Hilbert space. The Yukawa coupling matrix $D_F$ of the Standard Model provides the inverse line element for the finite geometry $(\cA_F,\cH_F,D_F)$ which displays the fine structure of space-time detected by the particles and forces discovered so far. For a long time the structure of the finite geometry was introduced ``by hand" following the trip inward bound \cite{Pais} in our understanding of matter and forces as foreseen by Newton: 
{\em
\begin{quote}
The Attractions of 	Gravity, Magnetism and Electricity, reach to very sensible distances, and so have been observed by vulgar Eyes, but there may be others which reach to so small distances as hitherto escape Observation.
\end{quote}}
and adapting by a bottom up process the finite geometry to the particles and forces, with the perfect fitting of the Higgs phenomenon and the see-saw mechanism with the geometric interpretation. There was however no sign of an ending in this quest, nor any sensible justification for the presence of the noncommutative finite structure. This state of affairs changed recently \cite{acmu1,acmu2} with the simultaneous  quantization of the fundamental class in $K$-homology and in $K$-theory. The $K$-homology fundamental class is represented by the Dirac operator. Representing the $K$-theory fundamental class, by requiring the use of the Feynman slash of the coordinates, explains the slight amount of non-commutativity of the finite algebra $\cA_F$ from Clifford algebras. From a purely geometric problem emerged the very same  finite algebra $\cA_F=M_2(\H) \oplus M_4(\C)$ which was the outcome of the bottom-up approach. 
 
One big plus of the noncommutative presentation is the economy in encoding. This economy is familiar in the written language 
where the order of letters is so important. Apart from the finite algebra $\cA_F$ the algebra of coordinates is generated by adjoining a ``punctuation symbol" $Y$ of $K$-theoretic nature, which is a unitary with $Y^4=1$. It is the noncommutativity of the obtained system which gets one outside the finite dimensional algebraic framework and generates the continuum. It allows one to write a higher analogue of the commutation relations which quantizes the volume and encodes all 4-dimensional spin geometries.

\end{mdframed}
\newpage
 \begin{mdframed}[backgroundcolor=blue!3]
 \centerline{\underline{\bf Zeta and NCG}}

\vspace{0.5cm} 
Let us  explain in a simple manner why the problem of locating the zeros of the Riemann zeta function is  intimately related to the ad\`ele class space. In fact what matters for the zeros is not the Riemann zeta function itself but rather the ideal it generates among holomorphic functions of a suitable class. A key role as a generator of this ideal is played by the operation on functions $f(u)$ of a real positive variable $u$ which is defined by
\begin{equation}\label{mapE}
	f \mapsto E(f), \ \ E(f)(v):=\sum_{\nt} f(nv)
\end{equation}
Indeed this operation is, in the variable $\log v$, a sum of translations $\log v\mapsto \log v+\log n $ by $\log n$  and thus after a suitable Fourier transform it becomes a product by the Fourier transform of the sum of the Dirac masses $\delta_{\log n}$,   \ie by the Riemann zeta function $\sum e^{-is\log n}=\zeta(is)$. It is thus by no means mysterious that one obtains a spectral realization by considering the cokernel of the map $E$. But what is interesting is the geometric meaning of this fact and the link with the explicit formulas. Thanks to the use of ad\`eles in the thesis of Tate the above summation over $\nt$ relates to a summation over  the  group $\Q^*$ of non-zero rational numbers. The natural adelic space on which the group $\Q^*$ is acting is simply the space $\A_\Q$ of ad\`eles over $\Q$ and in order to focus on the Riemann zeta function one needs to restrict oneself to functions invariant under the maximal compact subgroup $\hatz$ of the group of id\`ele classes which acts canonically on the noncommutative space of ad\`ele classes: $\Q^*\backslash\A_\Q$. In fact let $Y$ be the quotient $\A_f/\hatz$ of the locally compact space of finite ad\`eles by $\hatz$ and $y\in Y$ the element whose all components are $1$. The closure $F$ of the orbit $\nt y$ is compact and one has for $q\in \Q^*$ the equivalence
$
qy \in F\iff q\in \pm \nt
$.
 It is this property which allows one to replace the sum over $\nt$ involved in  \eqref{mapE} by the summation over the associated group $\Q^*$, simply by considering the function $1_F\otimes f$ on the product $Y\times \R$. The space of ad\`ele classes $X=\Q^*\backslash\A_\Q/\hatz=\Q^*\backslash(Y\times \R)$ is noncommutative because the action of $\Q^*$ on $\A_\Q$ is ergodic for the Haar measure. This geometrization gives a trace formula interpretation of the explicit formulas of Riemann-Weil (see \S \ref{sectzeta} below).

 The same space $X$ appears naturally from a totally different perspective which is that of Grothendieck toposes. Indeed the above action of $\nt$ by multiplication on the half line $[0,\infty)$ naturally gives rise to a 
Grothendieck topos and it turns out that $X$ is the space of points of this topos as shown in our recent work with C. Consani, \cite{CCscal,CCscal1}. Moreover the trace formula interpretation of the explicit formulas of Riemann-Weil allows one to rewrite the Riemann zeta function in the Hasse-Weil form where the role of the frobenius is played by its analogue in characteristic one. The space $X$ is the space of points of  the arithmetic site defined over the semifield $\rmax$, the Galois group of $\rmax$ over the Boolean semifield $\B$ is the multiplicative group $\R_+^*$ which acts by the automorphisms $x\mapsto x^\lambda$ and its action on $X$ replaces the action of the frobenius on the points over $\bar\F_q$ in finite characteristic. This topos theoretic interpretation yields a natural structure sheaf of tropical nature on the above scaling site  which becomes a curve in a tropical sense and the two sides : (topos theoretic, tropical) on one hand and (adelic, noncommutative) on the other hand, are in the same relative position as the two sides of the class field theory isomorphism.
\end{mdframed}
\vspace{0.2in}
The theory of Grothendieck toposes  is another crucial extension of the notion of space which plays a major role in algebraic geometry and is liable to treat quotient spaces in a successful manner.
The main difference of point of view between noncommutative geometry and the theory of toposes is that noncommutative geometry makes fundamental use of the complex numbers as coefficients and thus has a very close relation with the formalism of quantum mechanics. It takes complex numbers as the preferred coefficient field and Hilbert space as the main stage. The birth of noncommutative geometry can be traced back to this  night of June 1925 when, around three in the morning, W. Heisenberg while working alone in the Island of Helgoland in the north sea, discovered matrix mechanics and the noncommutativity of the phase space of microscopic mechanical systems. After the formulation of Heisenberg's discovery as matrix mechanics, von Neumann reformulated quantum mechanics using Hilbert space operators and went much further with Murray in identifying ``subsystems" of a quantum system as ``factorizations" of the underlying Hilbert space $\cH$. They discovered unexpected factorizations which did not correspond to tensor product decompositions $\cH=\cH_1\otimes \cH_2$ and developed the theory of factors which they classified into three types. In my thesis I showed that factors $\cM$ admit a canonical time evolution, \ie a canonical homomorphism 
$$
\R \stackrel{\delta}{\longrightarrow} {\rm Out}(\cM)=\Aut(\cM)/{\rm Int}(\cM)
$$
by showing that the class of the modular automorphism $\sigma_t^\varphi$ in ${\rm Out}(\cM)$ does not depend on the choice of the faithful normal state $\varphi$. This of course relied crucially on the Tomita-Takesaki theory \cite{TT} which constructs the map $\varphi \mapsto \sigma_t^\varphi$. The above uniqueness of the class of the modular automorphism \cite{Co_2bis}  drastically changed the status of the two invariants which I had previously introduced  in \cite{CoCR72a}, \cite{CoCR72b} by making them computable. The kernel of $\delta$,  
 $T(M)=\Ker \delta$ forms a subgroup of $\R$, the {\em periods} of $M$, and many non-trivial non-closed subgroups appear in this way. The fundamental invariant of factors is the {\em modular spectrum}  $S(M)$ of \cite{CoCR72a}. Its intersection $S(M)\cap \R_+^*$  is a closed subgroup of $\R_+^*$, \cite{acvd}, and this gave the subdivision of type III into 
 type  III$_\lambda$ $\iff$  $S(M)\cap \R_+^*=\lambda^Z$.  I showed moreover that the classification of factors of type  III$_\lambda$, $\lambda\in [0,1)$  is reduced to that of type II and automorphisms $$
M=N\rtimes_\theta \Z, \ \ N \ \textbf{type} \ {\rm II}_\infty, \ \theta \in \Aut(N)
$$
In the case III$_\lambda$, $\lambda\in (0,1)$, $N$ is a factor and the automorphism $\theta \in \Aut(N)$ is of module $\lambda$ \ie it scales the trace by the factor $\lambda$.
In the III$_0$ case, $N$ has a non-trivial center and, using the restriction of $\theta$ to the center, this gave a very rich invariant, a flow, and I used it in 1972, \cite{Cocras},  to show the existence of hyperfinite non ITPFI factor.  
 Only the case III$_1$ remained open in my thesis \cite{thesis} and was solved later by Takesaki \cite{Tak} using crossed product by $\R_+^*$.
 
 At that point one could fear that the theory of factors would remain foreign from main stream mathematics. This is not the case thanks to the following construction \cite{Co_foliation} of the von Neumann algebra $W(V,F)$ canonically associated to a foliated manifold $V$ with foliation $F$. It is the von Neumann algebra of random operators, \ie of bounded families $(T_\ell)$ parameterized by leaves $\ell$ and where each $(T_\ell)$ acts in the Hilbert space $L^2(\ell)$ of square integrable half-densities on the leaf. Modulo equality almost everywhere one gets a von Neumann algebra for the algebraic operations given by pointwise addition and multiplication: $$
W(V,F):=\{(T_\ell)\mid \text{acting on}\ L^2(\ell)\}
$$

 The  geometric examples lead  naturally to the most exotic factors, the prototype being the Anosov foliation of the sphere bundle of Riemann surfaces (see \S \ref{sectguillemin}) whose von Neumann algebra is the hyperfinite factor of type III$_1$.

\begin{figure}[H]
\begin{center}
\includegraphics[scale=0.45]{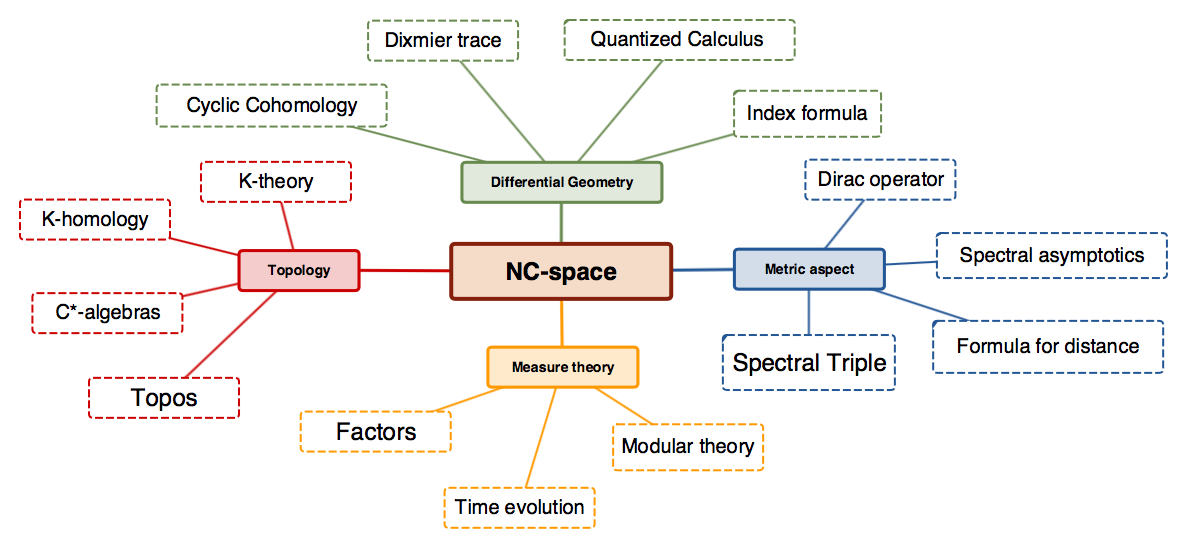}
\end{center}
\caption{Aspects of a noncommutative space \label{mindmap} }
\end{figure}

The above construction of the von Neumann algebra of a foliation only captures the ``measure theory" of the space of leaves. But such spaces inherit a much richer structure from the topological and differential geometric structure of the ambient manifold $V$. This fact served as a fundamental motivation for the development of noncommutative geometry and Figure \ref{mindmap} gives an overall picture of how the refined properties of leaf spaces and more general noncommutative spaces are treated in noncommutative geometry. One of the striking new features is the richness of the measure theory and how it interacts with a priori foreign features such as the secondary characteristic classes of foliations. Cyclic cohomology plays a central role in noncommutative geometry as the replacement of de-Rham cohomology \cite{ihes,pubihes}. One of its early striking applications obtained in \cite{transclass} (following previous work of S. Hurder)  using as a key tool the time derivative of the cyclic cocycle which is the transverse fundamental class, is that the flow of weights\footnote{The fundamental invariant of \cite{thesis,[Co-T]}} of the von Neumann algebra of a foliation with non-vanishing Godbillon-Vey class, admits an invariant probability measure. This implies immediately that the von Neumann algebra is of type III and exhibits the deep interplay between characteristic classes and ergodic theory of foliations.

\newpage
\section{New paradigm: spectral geometry}\label{sectspecgeom}

In the name ``Noncommutative Geometry" the adjective ``noncommutative" has an obvious algebraic meaning if one recalls that ``commutative" is one of the leitmotivs of algebraic geometry. As mentioned earlier reducing the new theory to the extension of previously known notions beyond the commutative case misses the point. In fact we shall now justify the term  ``Geometry" given to the field. The new paradigm for a geometric space is spectral. A spectral geometry is given by a triple $(\cA,\cH,D)$ where $\cA$ is an algebra of operators in the Hilbert space $\cH$ and where $D$ is a self-adjoint operator, generally unbounded, acting in the same Hilbert space $\cH$. Intuitively the algebra $\cA$ represents the functions on the geometric space and the operator $D$ plays the role of the inverse line element. It specifies how to measure distances by the formula 
\begin{equation}\label{dirac distance}
d(a,b)=\,{\rm Sup}\,\vert f(a)-f(b)\vert, f \ \text{such} \; \text{that}\ \Vert [D,f]\Vert\leq 1.
\end{equation}
Given a spin compact Riemannian manifold $M$ the associated spectral triple  $(\cA,\cH,D)$ is given by the algebra $\cA$ of functions on $M$ acting in the Hilbert space $\cH$ of $L^2$ spinors and the Dirac operator $D$. The formula \eqref{dirac distance} gives back the geodesic distance in the usual Riemannian sense but measures distances in a Kantorovich dual manner. The word ``spectral" in ``spectral triple" (or ``spectral geometry") has two origins, the straightforward one is the well-known Gelfand duality between a (locally compact) space and its algebra of complex valued continuous functions (vanishing at infinity). The deeper one is the full reconstruction of the geometry from the spectrum of  the Dirac operator (the $D$ in spectral triples) and the relative position in Hilbert space (the $\cH$ of the spectral triple) of the two von Neumann algebras given on one hand by the weak closure (double commutant) of the algebra $\cA$ (the $\cA$ of the spectral triple) and on the other hand the functions $f(D)$ of the operator $D$. We refer to \cite{ckm} for this point\footnote{One recovers the algebra of smooth functions from the von Neumann algebra of measurable ones  as those measurable functions which preserve the intersection of domains of powers of $D$. One then gets the topological space as the Gelfand spectrum of the norm closure of the algebra of smooth functions, one also gets its smooth structure, its Riemannian metric etc..}.  The great advantage of the new paradigm is that \eqref{dirac distance} continues to make sense for spaces which are no longer arcwise connected and works equally well for discrete or fractal spaces. Moreover the paradigm does not make use of the commutativity of $\cA$. When $\cA$ is noncommutative the evaluation $f(a)$ of $f\in\cA$ at a point $a$ no longer makes sense but the formula \eqref{dirac distance} still makes sense and measures the distance between states\footnote{The pure states are the extreme points of the space of states on $\cA$. In quantum mechanics they play the role of the points of the classical phase space. They are endowed with the natural equivalence relation defined by the unitary equivalence of the associated irreducible representations. The space of equivalence classes of pure states is typically a non-commutative space except in the easy type I situation.} on $\cA$, \ie positive linear forms $\phi:\cA\to \C$, normalized by $\phi(1)=1$. The formula measures the distance between states as follows :
\begin{equation}\label{dirac distance1}
d(\phi,\psi)=\,{\rm Sup}\,\vert \phi(f)-\psi(f)\vert, f \ \text{such} \; \text{that}\ \Vert [D,f]\Vert\leq 1.
\end{equation}
The new spectral paradigm of geometry extends the Riemannian paradigm and 
 we explained in \cite{CK} to which extent it allows one to go further in the query of Riemann on the  geometric paradigm of space in the infinitely small. Here ``space" is the space in which we live and is, after all, the one which truly deserves the name of ``geometry" independently of all the later developments which somehow usurp this name.  The new paradigm ensures that the line element does  encapsulate all the binding forces known so far
and is flexible enough to follow the inward bound quest of the geometric structure of space  for microscopic distances \cite{Pais}. The words of Riemann  are so carefully chosen and his vision so far-seeing that the reader will hopefully excuse our repeated use of this quotation:

{\em \begin{quotation}
Nun scheinen aber die empirischen Begriffe, in
welchen die r\"{a}umlichen Massbestimmungen gegr\"{u}ndet sind,
der Begriff des festen K\"{o}rpers und des Lichtstrahls, im
Unendlichkleinen ihre G\"{u}ltigkeit zu verlieren; es ist also
sehr wohl denkbar, dass die Massverh\"{a}ltnisse des Raumes im
Unendlichkleinen den Voraussetzungen der Geometrie nicht
gem\"{a}ss sind, und dies w\"{u}rde man in der That annehmen
m\"{u}ssen, sobald sich dadurch die Erscheinungen auf einfachere
Weise erkl\"{a}ren liessen. Es muss also entweder das dem Raume zu
Grunde liegende Wirkliche eine discrete Mannigfaltigkeit bilden,
oder der Grund der Massverh\"{a}ltnisse ausserhalb, in darauf
wirkenden bindenen Kr\"{a}ften, gesucht werden.\footnote{Now it seems that the
empirical notions on which the metric determinations of space
are founded, the notion of a solid body and of a ray of light,
cease to be valid in the infinitely small.  It is therefore
quite conceivable that the metric relations of space in
the infinitely small do not conform to the hypotheses of
geometry; and we ought in fact to assume this, if we can thereby
obtain a simpler explanation of phenomena. Either therefore the reality which
underlies space must be discrete, or we must
seek the foundation of its metric relations outside it, in binding
forces which act upon it.}
\end{quotation}}

We also explained in \cite{CK} the deep roots of the notion of spectral geometry in pure mathematics from the conceptual  understanding of the notion of (smooth, compact, oriented and spin) manifold. The fundamental property of a ``manifold" is not Poincar\' e duality in ordinary homology but is Poincar\' e duality in the finer theory called $KO$-homology. The key  result behind this, due to D. Sullivan (see \cite{MS}, epilogue), is that a PL-bundle is the same thing \footnote{Modulo the usual ``small-print" qualifications at the prime $2$, \cite{Siegel}} as a spherical fibration together with a $KO$-orientation.
The fundamental class in $KO$-homology contains all the information about the Pontrjagin classes of the manifold and these are not at all determined by its homotopy type: in the simply connected case only the signature class is fixed by the homotopy type.

  The second key step towards the notion of spectral geometry is the link with the Hilbert space and operator formalism which came as a byproduct of the work of Atiyah and Singer on the index theorem. They understood that operators in Hilbert space  provide the right realization for $KO$-homology cycles \cite{Atiyah, Singer}. Their original idea was developed by Brown-Douglas-Fillmore \cite{BDF}, Voiculescu \cite{Vo}, Mishchenko \cite{Mish}  and 
acquired its definitive form in the work of Kasparov at the end of the 1970's. The great new tool is  bivariant Kasparov theory \cite{Kas1,Kas2}, and  as far as  $K$-homology cycles are concerned the right notion is already in Atiyah's paper \cite{Atiyah}: A $K$-homology cycle on a compact space $X$ is given by a representation of the algebra $C(X)$ (of continuous functions on $X$) in a Hilbert space $\cH$, together with a Fredholm operator $F$ acting in the same Hilbert space fulfilling some simple compatibility condition (of commutation modulo compact operators) with the action of $C(X)$. One striking feature of this representation of $K$-homology cycles is that the definition does not make any use of the commutativity of the algebra $C(X)$.

 At the beginning of the 1980's, motivated by numerous examples of noncommutative spaces arising naturally in geometry from foliations and in particular the noncommutative torus $\T^2_\theta$ described in \S \ref{sectnctorus}, I realized that specifying an unbounded representative of the Fredholm operator gave the right framework for spectral geometry. The corresponding $K$-homology cycle\footnote{In \cite{BJ}, S. Baaj and P. Julg have extended the unbounded construction to the bivariant case. The most efficient axioms have been found by M. Hilsum in \cite{Hilsum} which allows to handle the case of symmetric non self-adjoint operators.} only retains the stable information and is insensitive  to deformations while the unbounded representative encodes the metric aspect. These are the deep mathematical reasons which are the roots of the notion of spectral triple. Since the year 2000  (\!\cite{ncg2000}) one important progress was obtained in the reconstruction theorem \cite{Crec} which gives an abstract characterization of the spectral triples associated to ordinary spin geometries.  Moreover a whole subject has developed concerning noncommutative manifolds and for lack of space we shall not cover it here but refer to the papers \cite{CDV1,CDV2,CDV3, CLa,DLSS,DVL1,DVL2,Landi3}. We apologise for the brevity of our account of the topics covered in the present survey and the reader will find a more complete treatment and bibliography in the collection of surveys \cite{chamssuijl,CCsurvey,DS5,Dundas, GE,FK3,GJV,LM1,LSZ1,Vo1,Wsurvey,XY2} to be published by Springer.

\subsection{Geometry and the modular theory} \label{sectnctorus}

\begin{figure}
\begin{center}
\includegraphics[scale=0.22]{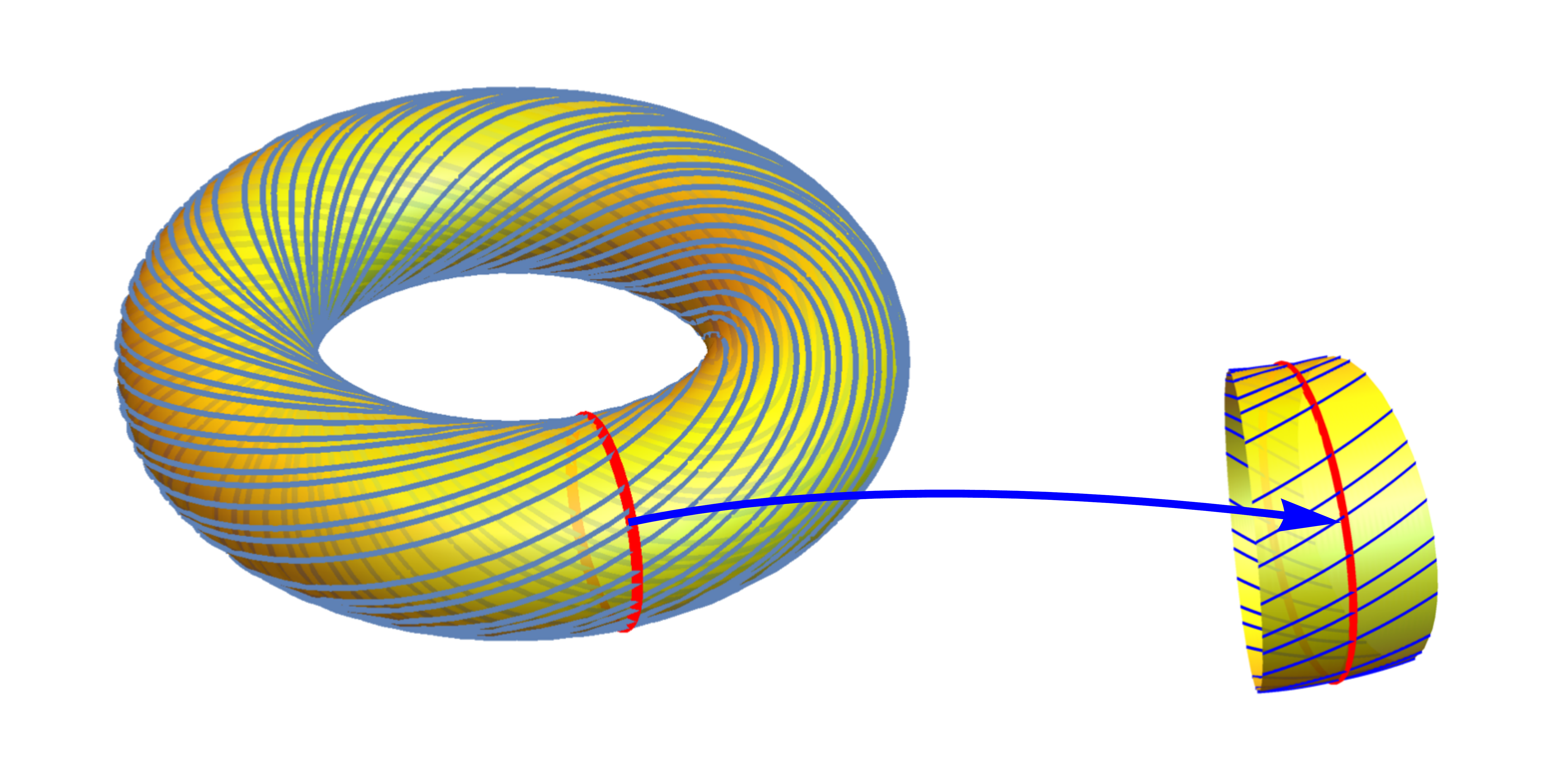}
\end{center}
\caption{Kronecker foliation $dy=\theta dx$ \label{fig:f1} }
\end{figure}

Foliated manifolds provide a very rich source of examples of noncommutative spaces and one is quite far from a full understanding of these spaces. In this section we discuss in an example the interaction of the new geometric paradigm with the modular theory (type III). The example is familiar from \cite{ncg2000}, it is  the noncommutative torus $\T^2_\theta$ whose $K$-theory was computed in the breakthrough paper of Pimsner and Voiculescu \cite{PV}. We first briefly recall the geometric picture underlying the noncommutative torus $\T^2_\theta$ and then describe the new development since year 2000 which displays a highly non-trivial interaction between the new paradigm of geometry and the modular type III theory. The underlying geometric picture is the Kronecker foliation $dy=\theta dx$ of the unit torus $V=\R^2/\Z^2$ as shown in Figure \ref{fig:f1}.  If one considers the restriction of the foliation to a neighborhood of a transversal as shown in Figure \ref{fig:f1} the structure is that of the product of a circle $S^1$ by an interval and the leaf space is $S^1$.  But for the whole foliation of the torus and  irrational $\theta$  the leaf space is a non-commutative space. It corresponds to the identification of those points on $S^1$ which are on the same orbit of the irrational rotation of angle $\theta$. At the algebraic level the foliation algebra is the same\footnote{up to Morita equivalence,  the algebraic  notion  has been adapted to $C^*$-algebras by M. Rieffel in\cite{Rie}}  as the crossed product $C(\T^2_\theta):=C(S^1)\rtimes_\theta \Z$  of the transversal  by the irrational rotation of angle $\theta$.

The differential geometry of the noncommutative torus $\T^2_\theta$ was first defined and investigated in  C.R. Acad. Sc. Paris Note of 1980 (\cf \cite{C}) ``$C^*$ alg\`ebres et g\'eom\'etrie diff\'erentielle''. 
While at the $C^*$-algebra level the presentation of the algebra $C(\T^2_\theta)$
of continuous functions on $\T^2_\theta$ is by means of two unitaries $U,V$ such that \[
VU = e^{2 \pi i \theta} \,U  V, 
\]
the smooth functions on $\T^2_\theta$ are given by the  dense subalgebra  $C^\infty(\T^2_\theta)\subset C(\T^2_\theta)$  of all series
\[
x =\sum_{m, n \in \mathbb{Z}} a_{m, n}\, U^m V^n, 
\]
such that the sequence of 
complex coefficients $(a_{m, n})$ is rapidly decaying in the sense 
that 
\[
\sup_{m, n \in Z} |a_{m, n}|(1 + |n| + |m|)^k < \infty, 
\]
for any non-negative integer $k$. Thus one has a simple explicit description of the generic element of $C^\infty(\T^2_\theta)$ (while there is no such explicit description for $C(\T^2_\theta)$). At the geometric level one has the smooth groupoid $G=S^1\rtimes_\theta \Z$ which is the reduction of the holonomy groupoid of the foliation to the transversal. The elements of $G$ are described equivalently as pairs of points of $S^1$ which are on the same orbit of the irrational rotation or as pairs $(u,n)\in S^1\times \Z$. The algebra is the convolution algebra of $G$. Any other transversal $L$ gives canonically a module $C^\infty(\T^2_\theta,L)$ over the convolution algebra of $G$ and the main discovery of \cite{C} can be summarized as follows:
\begin{enumerate}
\item The Schwartz space $\cS(\R)$ is a finite projective module over $C^\infty(\T^2_\theta)$ for the action of the generators $U,V$ on Schwartz functions $\xi(s)$ is given by translation of the variable $s\mapsto s-\theta$ and multiplication by $e^{2\pi i s}$.
\item The Murray-von Neumann dimension of $\cS(\R)$ is $\theta$.
\item The module	$\cS(\R)$ is naturally endowed with a connection of constant curvature and the product of the constant curvature by the Murray-von Neumann dimension is an integer.
\end{enumerate}

It is this integrality of the product of the constant curvature by the (irrational) Murray-von Neumann dimension which was the starting point of noncommutative differential geometry (\cf \cite{ihes,pubihes}). It holds irrespective of the transversal and the description of the general finite projective modules thus obtained was done in \cite{C}\footnote{and was often misattributed as in \cite{manin}}.
 The differential geometry as well as pseudo-differential operator calculus of the noncommutative two torus $\T^2_\theta$ were first developed in \cite{C}. As in the case of an ordinary torus the conformal structure is specified by a complex modulus $\tau\in \C$ with $\Im(\tau)>0$. The flat spectral geometry $(\cA,\cH,D_0)$ is constructed using the two derivations $\delta_1,\delta_2$ which generate the action by rotation of the generators $U,V$ of the algebra $\cA=C^\infty(\T^2_\theta)$. To obtain a curved geometry $(\cA,\cH,D)$ from the flat one $(\cA,\cH,D_0)$,  one introduces (\cf \cite{cc}, \cite{Paula}) a noncommutative Weyl conformal factor (or dilaton), 
which changes the metric by modifying the noncommutative volume form 
while keeping the same conformal structure. 
Both notions of volume form and of conformal structure are well understood in the general case (\cf \cite[\S VI]{Co-book}). The new and
crucial ingredient which has no classical analogue is 
 the modular operator $\modu$ of the non-tracial weight $\vp(a)=\vp_0(a e^{-h})$ associated to the dilaton $h$.

\medskip

 In noncommutative geometry the local geometric invariants such as the Riemannian curvature are extracted from the spectral triple using 
 the functionals defined by the coefficients of heat kernel expansion 
$$
\Tr(a e^{-tD^2})\, \sim_{t \searrow 0} \, \sum_{n \geq 0} {\rm a}_n(a,D^2)t^{\frac{-d+n}{2}} \, , \quad a \in \cA ,
$$
where $d$ is the dimension of the geometry. Equivalently, 
one may consider special values of the corresponding zeta functions.
 Thus, the local curvature of the geometry is detected by the high frequency behavior of the spectrum of $D$ coupled with the action of the algebra $\cA$ in $\cH$. 
 In the case of $\T^2_\theta$ the  computation of the value at $s=0$ of the zeta function $\Tr(a|D|^{-2s})$ for the $2$-dimensional
curved geometry associated to the dilaton $h$, or equivalently of the coefficient
${\rm a}_2(a,D^2)$ of the heat expansion  was
started in the late 1980's (\cf \cite{cc}), and the proof of the analogue of the Gauss--Bonnet formula was published in \cite{Paula}. It was subsequently extended in \cite{FK} 
to the case of arbitrary values of the complex modulus $\tau$ (one had $\tau=i$ in
\cite{Paula}). In these papers only the total integral of the curvature was needed, and this allowed one to make simplifications under the trace which are no longer possible  when one wants to fully compute the local expression for the functional
$a \in \cA \mapsto {\rm a}_2(a,D^2)$. 

 \medskip

The complete calculation  of ${\rm a}_2(a,D^2)$ was actually performed in our joint work with H. Moscovici 
in \cite{ConMosModular}, and  
 independently  by F. Fathizadeh and M. Khalkhali in
\cite{FK1}. The spectral geometry of noncommutative tori is a key testground for  new ideas and basic contributions have been done in \cite{Paula,ConMosModular,LM,FK2,FGK,GKh, Khal,FK3,FGH,Liu,IM,Ponge3,Lesch}, with beautiful surveys by M. Lesch and H. Moscovici \cite{LM1} and by  F. Fathizadeh and M. Khalkhali \cite{FK3} to which we refer for a more complete account and bibliography. In this topic the hard direct calculations alternate with their conceptual understanding and the main result  of \cite{ConMosModular} is  the conceptual explanation for the complicated formulas. The explanation is obtained by expressing 
in terms of a closed formula the 
Ray-Singer log-determinant of $D^2$. 
The gradient of the log-determinant functional yields in turn a local curvature formula,
which arises as a sum of two terms, each involving a 
function in the modular operator, of one and respectively two variables.
Computing the gradient in two different ways leads to the proof of
a deep internal consistency relation between these two distinct constituents. This relation was checked successfully and gave confidence in the pertinence of the whole framework.  Moreover it elucidates the meaning of the intricate two operator-variable function which we now briefly describe. As in the case of the standard torus viewed as a complex curve,
the total Laplacian associated to such a spectral triple splits into two
components, one $\perta$ on functions and the other $\pertc$  on $(0,1)$-forms, the two
operators being isospectral outside zero. The
corresponding curvature formulas involve
second order (outer) derivatives of the Weyl factor.
 For the Laplacian $\perta$ the result is of the form
\begin{equation}\label{a2term}
    {\rm a}_2(a,\perta)=-\frac{\pi}{2\tau_2}\vp_0(a\left(K_0(\nabla)(\lapa(h))+\frac 12 H_0(\nabla_1,\nabla_2)(\square_\Re (h)\right) ,
\end{equation}
where $\nabla=\log \modu$ is the inner derivation implemented  by $-h$,
$$
 \lapa(h)=
  \delta_1^2(h)+2 \Re(\tau)\delta_1\delta_2(h)+|\tau|^2 \delta_2^2(h) ,
 $$
 $\square_\Re$ is the Dirichlet quadratic form
 $$
 \square_\Re (\ell) :=
(\delta_1(\ell))^2+ \Re(\tau)\left(\delta_1(\ell)\delta_2(\ell)+\delta_2(\ell)\delta_1(\ell)\right)+|\tau|^2 (\delta_2(\ell))^2\,,
$$
and $\nabla_i, \, i=1, 2$, signifies that $\nabla$ is acting on the $i$th factor.
The operators $K_0(\nabla)$ and $H_0(\nabla_1,\nabla_2)$ are new ingredients, whose
 occurrence has no classical analogue and is a vivid manifestation of the  non-unimodular nature of
 the  curved geometry of the noncommutative $2$-torus with Weyl factor.
 The functions $K_0(u)$ and $H_0(u,v)$ by which the modular derivatives act
  seem at first very complicated, and of course beg for a
  conceptual understanding. This was obtained in \cite{ConMosModular} where  denoting
$
\tilde{K}_0(s)\, = \, 4\frac{\sinh(s/2)}{s} K_0(s)$ and 
$\tilde{H}_0(s,t)\, = \, 4\frac{\sinh((s+t)/2)}{s+t} H_0(s,t),
$ 
we give an abstract proof of the 
 functional relation  
\begin{equation}\label{transforulepre}
  - \frac 12 \tilde H_0(s_1,s_2)=\frac{\tilde K_0(s_2)-\tilde K_0(s_1)}{s_1+s_2}+\frac{\tilde K_0(s_1+s_2)-\tilde K_0(s_2)}{s_1}-\frac{\tilde K_0(s_1+s_2)-\tilde K_0(s_1)}{s_2}
\end{equation}
which determines the whole structure in terms of the function $\tilde K_0$ which turns out to be (up to the factor $\frac 18$) the generating function of the Bernoulli numbers, \ie one has $\frac 18\tilde K_0(u)=\sum_1^\infty \frac{B_{2n}}{(2n)!}u^{2n-2}\,.$

Among the recent developments are the full understanding of the interplay of the geometry with Morita equivalence (a purely noncommutative feature) in \cite{LM} and a new hard computation in \cite{CFF}. This latter result gives the formulas for the $a_4$ term for the noncommutative torus and, as a consequence, allows one to compute this term for the noncommutative $4$-tori which are products of two $2$-tori. We managed in \cite{CFF} to derive abstractly and check the analogues of the above functional relations, but the formulas remain very complicated and mysterious, thus begging for a better understanding. One remarkable fact is that  the noncommutative $4$-tori which are products of two $2$-tori, are no longer conformally flat so that the computation of the $a_4$ term for them involves a more complicated modular structure which is given by a two parameter group of automorphisms. This provides a first hint in the following program:

\vspace{0.5cm}
 \begin{mdframed}[backgroundcolor=red!8]
The long range goal in exploring the spectral geometry of noncommutative tori is to arrive at the formulation of the full fledged modular theory involving the full Jacobian matrix rather than its determinant, as suggested long ago by the transverse geometry of foliations in \cite{transclass} where the reduction from type III to type II was done for the full jacobian. The hypoelliptic theory was used in \cite{CMos} to obtain the transverse geometry at the type II level. The first step of the general theory is to adapt the twist defined in \cite{ConMosTypeIII} to a notion involving the full Jacobian.
\end{mdframed}

\subsection{Inner fluctuations of the metric}\label{sectinnerfluct}

In the same way as inner automorphisms enhance the group of diffeomorphisms,  they provide special deformations of the metric for a spectral geometry $(\cA,\cH,D)$ which correspond to the gauge fields of non-abelian gauge theories. 
In our joint work with A. Chamseddine and W. van Suijlekom \cite{acinner}, we obtained a conceptual understanding of the role of the gauge bosons in physics as the inner fluctuations of the metric. I will describe this result here in a non-technical manner. 

 Ignoring  the important nuance coming from the real structure $J$, the inner fluctuations of the metric were first defined as the transformation 
$$
D\mapsto D+A, \ \ A=\sum a_j[D,b_j], \ \ a_j,  b_j \in \cA, \ A=A^*
$$ 
which imitates the way classical gauge bosons appear as matrix-valued one-forms in the usual framework. The really important facts were that the spectral action applied to $D+A$ delivers the Einstein-Yang-Mills action which combines gravity with matter in a natural manner, and that the gauge invariance becomes transparent at this level since an inner fluctuation coming from a gauge potential of the form $A=u[D,u^*]$ where $u$ is a unitary element (\ie $uu^*=u^*u=1$) simply results in a unitary conjugation $D\mapsto uDu^*$ which does not change the spectral action.
What we discovered in our joint work with A. Chamseddine and W. van Suijlekom \cite{acinner} is that the inner fluctuations arise in fact from the action on metrics (\ie the $D$) of a canonical {\em semigroup} $\pert(\cA)$ which only depends upon the algebra $\cA$ and extends the unitary group.  The semigroup is defined as the self-conjugate elements:
$$
\pert(\cA):=\{A=\sum a_j\otimes b_j^{\rm op}\in  {\mathcal{A}}\otimes {\mathcal{A}}^{\rm op}\mid \sum a_jb_j=1, \ \theta(A)=A\}
$$ 
where $\theta$ is the antilinear automorphism of the algebra ${\mathcal{A}}\otimes {\mathcal{A}}^{\rm op}$ given by 
$$
\theta:\sum a_j\otimes b_j^{\rm op}\mapsto \sum b_j^*\otimes a_j^{*\rm op}.
$$
The composition law in $\pert(\cA)$ is the product in the algebra ${\mathcal{A}}\otimes {\mathcal{A}}^{\rm op}$. The action of this semigroup $\pert(\cA)$ on the metrics is given, for $A=\sum a_j\otimes b_j^{\rm op}$  by 
$$
D\mapsto D'=^{A}\!\!D= \sum a_j D b_j.
$$
Moreover, the transitivity of inner fluctuations which is a key feature since one wants that an inner fluctuation of an inner fluctuation is itself an inner fluctuation now results from the semigroup structure according to the equality:
$
^{A'}\!(^{A}D)=^{(A'A)}\!D.
$
We refer to \cite{acinner} for the full treatment. 
Inner fluctuations have been adapted in the twisted case in \cite{Landi1,Landi2}.

\subsection{The spectral action and standard model coupled to gravity}

The 
 new spectral paradigm of geometry,  because of its flexibility, provides the needed tool to 
 refine our understanding of the structure of physical space in the small and to ``{\it seek the foundation of its metric relations outside it, in binding
forces which act upon it}". The main idea, described in details in \cite{CK}, is that the line element now embodies not only the force of gravity but all the known forces,  electroweak and strong, appear from the spectral action and the inner fluctuations of the metric. This provides   a completely new perspective on the geometric interpretation of the detailed  structure of the Standard model and of the Brout-Englert-Higgs mechanism. One gets the following simple mental picture for the appearance of the scalar field : imagine that the space under consideration is two sided like a sheet $S$ of paper in two dimensions. Then when differentiating a function on such a space one may restrict the function to either side $S_\pm$ of the sheet and thus obtain two spin one fields. But one may also take the finite difference $f(s_+)-f(s_-)$ of the function at the related points of the two sides. The corresponding field is clearly insensitive to local rotations and is a scalar spin zero field. This, in a nutshell, is how the Brout-Englert-Higgs  field appears geometrically once one accepts that there is a ``fine structure" which is revealed by the detailed structure of the standard model of matter and forces. This allows one to uncover  the  
geometric meaning of the Lagrangian of gravity coupled to the standard model. This extremely complicated Lagrangian is obtained from the spectral action developed in our joint work with A. Chamseddine \cite{cc2} which is the only natural additive spectral invariant of a noncommutative geometry. In order to comply with Riemann's requirement that the inverse line element $D$ embodies the forces of nature, it is evidently important that we do not separate artificially the gravitational part from the gauge part, and that $D$ encapsulates both forces in a unified manner. In the traditional geometrization of physics, the gravitational part specifies the metric while the gauge part corresponds to a connection on a principal bundle. In the NCG framework, $D$ encapsulates both forces in a unified manner and the gauge bosons appear as inner fluctuations of the metric but form an inseparable part of the latter. The noncommutative geometry dictated by physics is the product of the ordinary $4$-dimensional continuum by a finite noncommutative geometry $(\cA_F,\cH_F,D_F)$ which appears naturally from the classification of finite geometries of $KO$-dimension  equal to $6$ modulo $8$ (\cf \cite{cc5,mc2}). The finite dimensional algebra $\cA_F$ which appeared is of the form 
$$
\cA_F=C_{+}\oplus C_{-}, \ \  C_{+}=M_{2}(\mathbb{H}),\  \  \  C_{-}=M_{4}(\mathbb{C}).
$$

The agreement of the mathematical formalism of spectral geometry and all its subtleties such as the periodicity of period $8$  of the $KO$-theory, might still be accidental but I personally  got convinced of the pertinence of this approach when recovering \cite{mc2} the see-saw mechanism (which was dictated by the pure math calculation of the model) while I was unaware of its key physics role to provide very small non-zero masses to the neutrinos, and of how it is ``put by hand" in the standard model. The low Higgs mass then came in 2012 as a possible flaw of the model, but in \cite{acresil} we proved the compatibility of the model with the measured value of the Higgs mass,  due to the role in the renormalization of the scalar field which was already present in  \cite{ac2010} but had been ignored thinking that it would not affect the running of the self-coupling of the Higgs.

In all the previous developments we had followed the ``bottom-up" approach, \ie we uncovered the details of the finite noncommutative geometry $(\cA_F,\cH_F,D_F)$ from the experimental information contained in the standard model coupled to gravitation. In 2014, 
in collaboration with A. Chamseddine and S. Mukhanov \cite{acmu1,acmu2} we were investigating  the purely geometric problem of encoding $4$-manifolds in the most economical manner in the spectral formalism. Our problem had no a priori link with the standard model of particle and forces and the idea was to treat the coordinates in the same way as the momenta are assembled together in a single operator using the gamma matrices. The great surprise was that this investigation  gave the conceptual explanation of the finite noncommutative geometry from Clifford algebras! This is described in details in \cite{CK} to which we refer. What we  obtained is a higher form of the Heisenberg commutation relations between $p$ and $q$, whose irreducible Hilbert space representations correspond to $4$-dimensional spin geometries. The role of $p$ is played by the 
Dirac operator and the role of $q$ by the Feynman slash of coordinates using Clifford algebras. The proof that all spin geometries are obtained relies on deep results of immersion theory and ramified coverings of the sphere.  The volume of the $4$-dimensional geometry is automatically quantized by the index theorem and the  spectral model, taking into account the inner automorphisms due to the noncommutative nature of the Clifford algebras, gives Einstein gravity coupled with the slight extension of the standard model which is a Pati-Salam model. This model was shown in \cite{acpati1,acpati2} to yield unification of coupling constants. We refer to the survey \cite{chamssuijl} by A. Chamseddine and W. van Suijlekom for an excellent account of the whole story of the evolution of this theory from the early days to now. 

 The dictionary between the physics terminology and the fine structure of the geometry is of the following form:

\vspace{0.5cm}
 \begin{mdframed}[backgroundcolor=green!3]
\begin{center}
\begin{tabular}{r | l}

 {\bf Standard Model} & {\bf Spectral Model}\\
  & \\
 Higgs Boson &  Inner metric$^{(0,1)}$ \\
 &\\
   Gauge bosons  &  Inner metric$^{(1,0)}$ \\
  &\\
     Fermion masses $u,\nu$&     Dirac$^{(0,1)}$ in $\uparrow$\\
    &\\
    CKM matrix,  Masses down &   Dirac$^{(0,1)}$ in $\dodo$\\
    &  \\
  Lepton mixing, Masses leptons $e$ &   Dirac$^{(0,1)}$ in $\elel$\\
  & \\
         Majorana  mass matrix & Dirac$^{(0,1)}$ on $E_R\oplus J_F E_R$\\
   & \\
  Gauge couplings  & Fixed at unification\\
    &  \\
  Higgs scattering parameter   &  Fixed at unification  \\
     &  \\
  Tadpole constant  &  $- \mu_0^2\, |\higgs|^2$\\
   &  \\
  \end{tabular}
\end{center}
\end{mdframed}
\vspace{0.5cm}

\subsection{Dimension $4$}\label{sectsobolev}

In our work with A. Chamseddine and S. Mukhanov \cite{acmu1,acmu2}, the dimension $4$ plays a special role for the following reason. In order to encode a manifold $M$ one needs to construct a pair of maps $\phi,\psi$ from $M$  to the sphere (of the same dimension $d$) in such a way that the sum of the pullbacks
of the volume form of the (round) sphere vanishes nowhere on $M$. This problem is easy to solve in dimension $2$ and $3$ because one first writes $M$ as a ramified cover $\phi:M\to S^d$  of the sphere and one pre-composes   $\phi$ with a diffeomorphism $f$ of $M$ such that $\Sigma\cap f(\Sigma)=\emptyset$ where $\Sigma$ is the subset where $\phi$ is ramified. This subset is of codimension $2$ in $M$ and there is no difficulty to find $f$ because $(d-2)+(d-2)<d$ for $d<4$. It is worth mentioning that the $2$ for the codimension of $\Sigma$ is easy to understand from complex analysis: For an arbitrary smooth map $\phi:M\to S^d$ the Jacobian will vanish on a codimension $1$ subset, but in one dimensional  complex analysis the Jacobian is a sum of squares and its vanishing means the vanishing of the derivative which gives two conditions rather than one. Note also in this respect that quaternions do not help.\footnote{As an exercice one can compute the Jacobian of the power map $q\mapsto q^n$, $q\in \H$, and show that it vanishes on a codimension $1$ subset. For instance for $n=3$ the Jacobian is given by $9 \left(3 a^2-b^2-c^2-d^2\right)^2 \left(a^2+b^2+c^2+d^2\right)^2$.}  Thus dimension $d=4$ is the critical dimension for the above existence problem of the pair $\phi,\psi$. Such a pair does not always exist\footnote{It does not exist for $M=P_2(\C)$.} but as shown in \cite{acmu1,acmu2}, it always exist for spin manifolds which is the relevant case for us. 

The higher form of the Heisenberg commutation relation mentioned above involves in dimension $d$ the power $d$ of the commutator $[D,Z]$ of the Dirac operator with the operator $Z$ which is constructed (using the real structure $J$,  see \cite{CK}) from the coordinates. We shall now explain briefly how this fits perfectly with the framework of D. Sullivan on Sobolev manifolds, \ie of manifolds of dimension $d$ where the pseudo-group underlying the atlas  preserves continuous functions with one derivative in $L^d$.   He discovered  the intriguing special role of dimension $4$ in this respect. He showed in \cite{S3} that 
 topological manifolds in dimensions $>5$ admit bi-lipschitz coordinates and these are unique up to small perturbations, moreover existence and uniqueness also holds for  Sobolev structures: one derivative in $L^d$.
A stronger result was known classically for dimensions 1,2 and 3. There the topology controls the smooth structure up to small deformation. 
In dimension 4, he proved with S. Donaldson in \cite{DS} that for manifolds with coordinate atlases related by the  pseudo-group  preserving continuous functions with one derivative in $L^4$ it is possible to  develop the $SU(2)$ gauge theory 
and the  famous Donaldson invariants. Thus in dimension $4$ the Sobolev manifolds  behave like the smooth ones as opposed to Freedman's abundant topological manifolds.\footnote{for which  any modulus of continuity  whatsoever is not known.} The obvious question then is to which extent  the higher Heisenberg equation of \cite{acmu1,acmu2} singles out the Sobolev manifolds as the relevant ones for the functional integral involving the spectral action. 

\subsection{Relation with  quantum gravity} 
Our approach on the geometry of space-time is not concerned with quantum gravity but it addresses a more basic related question which is to understand ``why gravity coupled to the standard model", which we view as a preliminary. The starting point of our approach  is an extension of Riemannian geometry beyond its classical domain which provides the needed flexibility in order to answer the querry of Riemann in his inaugural lecture. The modification of the geometric paradigm comes from the confluence of the abstract understanding of the notion of manifold from its fundamental class in $KO$-homology with the advent of the formalism of quantum mechanics. This confluence comes from the realization of cycles in $KO$-homology from representations in Hilbert space. The final touch on the understanding of the geometric reason behind gravity coupled to the standard model, came from the simultaneous quantization of the fundamental class in $KO$-homology and its dual in $KO$-theory which gave rise to the  higher Heisenberg relation. So far we remain at the level of first quantization and the problem of second quantization is open. A reason why this issue cannot be ignored is that the quantum corrections to the line element as explained in Section \ref{sectlocalact} below (and shown in Figure \ref{dsdr}) are ony the tip of the iceberg since the dressing also occurs for all the $n$-point functions for Fermions while Figure \ref{dsdr} only involves  the two point function. One possible strategy to pass to this second quantized higher level of geometry is to try to give substance to the proposal we did in \cite{CK}:
{\it
\begin{quote}
The duality between $KO$-homology and $KO$-theory is the origin of the higher Heisenberg relation. As already mentioned in \cite{coinaugural}, algebraic $K$-theory, which is a vast refinement of topological $K$-theory,  is begging for the development of a dual theory and one should expect profound relations between this dual theory and the  theory of interacting quanta of geometry. As a concrete point of departure, note that the deepest results on the topology of diffeomorphism groups of manifolds are given by the Waldhausen algebraic $K$-theory of spaces and we refer to \cite{DGM} for a unifying picture of algebraic $K$-theory. 	
\end{quote}}
\vspace{0.5cm}
 \begin{mdframed}[backgroundcolor=red!8]
The role of the second quantization of Fermions in `` Entropy and the spectral action"   \cite{CCS} which interprets the spectral action as an entropy  is a first step towards building from quantum field theory a {\em second quantized} version of spectral geometry. 
\end{mdframed}
\vspace{0.5cm}

Note also that an analogous  result to \cite{CCS} has been obtained in \cite{DK} for the bosonic case.

\section{New tool: Quantized calculus}

A key tool of noncommutative geometry is the quantized calculus. Its origin is the notion of real variable provided by the formalism of quantum mechanics in terms of self-adjoint operators in Hilbert space. One basic very striking fact is that if $\phi:\R\to \R$ is an arbitrary Borel function and $H$ a self-adjoint operator in Hilbert space, then $\phi(H)$ makes sense even though $\phi$ can have a completely different definition in  various parts of $\R$. It is this fact which reveals in which sense  the formalism of quantum mechanics is compatible with the naive notion of real variable. It is superior to the classical notion of random variable because it models discrete variables with finite multiplicity as coexisting with continuous variables while a discrete random variable with finite multiplicity only exists on a countable sample space. The price one pays for this coexistence is non-commutativity since the commutant of a self-adjoint operator $H$ with countable spectrum of multiplicity $1$ is the algebra of Borel functions  $\phi(H)$ and they all have countable range and atomic spectral measure.  This implies in particular that the infinitesimal variables of the theory, which are modeled by compact operators, cannot commute with the continuous variables. We refer to \cite{CK} for a detailed discussion of this point. We now display the dictionary which translates the classical geometric notions into their noncommutative analogue.

\subsection{Quantum formalism and variables, the dictionary}

We refer to Chapter IV of \cite{Co-book} for a detailed description of the calculus and many concrete examples. 

\vspace{0.5cm}
 \begin{mdframed}[backgroundcolor=yellow!3]
\begin{tabular}{r | l}
& \\
 Real variable $f:X\to \R$  & Self-adjoint operator $H$ in Hilbert space  \\
Range $f(X)\subset \R$ of the variable & Spectrum of  the operator $H$\\ 
Composition $\phi\circ f$, $\phi$ measurable & Measurable functions $\phi(H)$ of self-adjoint operators  \\
Bounded complex variable $Z$ & Bounded operator $A$ in Hilbert space \\
Infinitesimal variable $dx$ & Compact operator $T$\\
Infinitesimal of order  $\alpha>0$  & Characteristic values $\mu_n(T)=O(n^{-\alpha})$ for $n\to \infty$ \\
 Algebraic operations on functions & Algebra of operators in Hilbert space \\ 
Integral of function $\int f(x)dx$ & $\displaystyle{\int\!\!\!\!\!\! -} T =$ coefficient of $\log(\Lambda)$ in $\Tr_\Lambda(T)$\\
Line element $ds^2=g_{\mu\nu} dx^\mu dx^\nu$ & $ds=\bullet\!\!\!\!-\!\!\!-\!\!\!-\!\!\!-\!\!\!-\!\!\bullet$ : Fermion propagator  $D^{-1}$\\
$d(a,b)=\,{\rm Inf}\,\int_\gamma\,\sqrt{g_{\mu\,\nu}\,dx^\mu\,dx^\nu}$ & $d(\mu,\nu)=\,{\rm Sup}\,\vert \mu(A)-\nu(A)\vert, \mid \ \Vert [D,A]\Vert\leq 1.$ \\ 
Riemannian geometry $(X,ds^2)$& Spectral geometry $(\cA,\cH,D)$\\
Curvature invariants & Asymptotic expansion of spectral action \\
Gauge theory & Inner fluctuations of the metric\\
Weyl factor perturbation & $D\mapsto \rho D\rho$ \\
Conformal Geometry & Fredholm module $(\cA,\cH,F)$, $F^2=1$. \\
Perturbation by Beltrami differential & $F\mapsto (aF+b)(bF+a)^{-1}$, $a=(1-\mu^2)^{-1/2}$, $b=\mu a$ \\
Distributional derivative & Quantized differential $dZ:=[F,Z]$ \\
Measure of conformal weight $p$ & $f \mapsto \displaystyle{\int\!\!\!\!\!\! -}f(Z)\vert dZ\vert^p$ \\
&
\end{tabular}
\end{mdframed}

\vspace{0.5cm}

While the first lines of the dictionary are standard in quantum mechanics, they become more and more involved as one goes down the list. We refer to the survey of Dan Voiculescu \cite{Vo1} for his deep work on commutants modulo normed ideals and the modulus of quasicentral approximation for an $n$-tuple of operators with respect to a normed ideal which provides a fundamental tool, used crucially in the reconstruction theorem \cite{Crec}. We did  include the inner fluctuations in the dictionary as corresponding to the gauge theories and the precise relation is explained in \S \ref{sectinnerfluct}. Note also in this respect that the alteration of the metric by a Weyl factor is of a similar nature as the inner fluctuations. One understands from \cite{Co-book}, how the Beltrami differentials act on conformal structures when the latter are encoded in bounded form, but the problem of encoding a general change of metric is yet open. After listing below the various entries of the dictionary we shall concentrate on the highly non-trivial part which is integration in the form of the Dixmier trace and the Wodzicki residue. 

\subsection{The principle of locality in NCG}\label{sectlocal}
 It is worth describing a fundamental principle which has emerged over the years. It concerns the notion of ``locality" in noncommutative geometry. The notion of locality is built in for  topological spaces or more generally for toposes. It is also an essential notion in modern physics. 
 
 Locality also plays a key role in noncommutative geometry and makes sense in the above framework of the quantized calculus, but it acquires another  meaning, more subtle than the straightforward topological one. The idea 
 comes from the way the Fourier transform translates the local behavior of functions on a space in terms of the decay at $\infty$ of their Fourier coefficients.  The relevant properties of the ``local functionals" which makes them ``local" is translated in Fourier by their dependance on the asymptotic behavior of the Fourier coefficients.
 
  This is familiar in particle physics where it is the ultraviolet behavior which betrays the fine local features, in particular through the ultraviolet divergencies. More specifically, in Quantum Field Theory, the counterterms
 which need to be added along the way in the renormalization process  as corrections to the initial Lagrangian appear from the ultraviolet divergencies and are automatically ``local".  One merit of the divergencies is to be insensitive  to bounded perturbations of the momentum space behavior.   
 The fundamental ones appear either as the coefficient of a logarithmic divergency when a cutoff parameter $\Lambda$ (with the dimension of an energy) is moved to $\infty$ or, in the Dim-Reg regularization scheme, as a residue. 
 
 In noncommutative geometry the above two ways of isolating the relevant term in divergencies have been turned into two key tools which allow one to construct ``local functionals". These are 
 
 \vspace{0.2in}
$\blacktriangleright$~The Dixmier trace

\vspace{0.2in}
$\blacktriangleright$~The Wodzicki residue
 
 \vspace{0.2in}
 
 These tools are functionals on suitable classes of operators, and their role is to filter out the irrelevant details. In some sense these filters  wipe out the quantum details and give a semi-classical picture.

\subsubsection{Dixmier trace and fractals}

The integration in noncommutative geometry needs to combine together two properties which seem contradictory: being a positive trace (in order to allow for cyclic permutations under the functional) and vanishing on infinitesimals of high enough order (order $>1$). This latter property together with the approximation from below of positive compact operators by finite rank ones entails that the integration cannot commute with increasing limits and is non-normal as a functional. Such a functional was discovered by J. Dixmier in 1966, \cite{Di}. One might be skeptical at first, on the use of such functionals for the reason which was invoked in the introduction concerning the non-constructive aspect of some proofs for the cardinality of natural non-commutative spaces. Fortunately there is a perfect answer to this objection both at the abstract general level as well as in practice. The ingredient needed to construct the non-normal functional is a limiting process $\lim_\omega$ on bounded sequences of real numbers. It is required to be linear and positive (\ie  $\lim_\omega((a_n))\geq 0$ if $a_n\geq 0$ $\forall n$) and to vanish for sequences which tend to $0$ at infinity. The book of Hardy \cite{Har} on divergent series shows the role of iterated Cesaro means and of such limiting processes in many questions of number theory related to Tauberian theorems. If one would require the limiting process to be multiplicative it would correspond to an ultrafilter and would easily be shown to break the rule of measurability which we imposed from the start in the introduction. At this point a truly remarkable result of G. Mokobodzki \cite{Mey} saves the day. He managed to show, using the continuum hypothesis, the existence of a limiting process which besides the above conditions is {\bf universally measurable} and fulfills 
\begin{equation}\label{dixmeasure}
\lim_{\omega} \left( \int a(\alpha)d\mu(\alpha)\right ) =\int  \lim_{\omega}(a(\alpha))d\mu(\alpha)
\end{equation}
for any bounded measurable family $a(\alpha)$ of sequences of real numbers. In other words if one uses this medial limit of Mokobodski one no longer needs to worry about the measurability problems of the limit and one can permute it freely with integrals! One might object that the medial limit is not explicitly constructed since its existence relies on the continuum hypothesis but the medial limit is to be used as a tool and there are very few places in mathematics where one can use the power of subtle independent axioms such as the continuum hypothesis. We refer to the article of Paul Cohen \cite{PC} for the interplay between axioms of set theory and mathematical proof. He explains how certain proofs might require the use of stronger axioms even though the original question is formulated in simple number theoretic language. But the objection that the Dixmier trace is ``non-constructively defined" does not apply in practice, the main point there being that this issue does not arise when there is convergence. Since year 2000 a lot of progess has been done and for instance the  construction of Dixmier traces and the proof of its important properties have been extended to the type II situation in \cite{BF}. Moreover a very powerful school with great analytical skill has developed around F. Sukochev and D. Zanin and their book with S. Lord \cite{LSZ}. With them we have undertaken the task of giving complete proofs of several important results of \cite{Co-book} which were only briefly sketched there.  In \cite{CSZ}  we give complete proofs of Theorem 17 of Chapter IV of \cite{Co-book} and in \cite{CSZ1}, with E. McDonald, the proof of the result on Julia sets announced in  \cite{Co-book} page 23. At the technical level an important hypothesis plays a role concerning the invariance of the limiting process under rescaling. It is simpler to state it 
 for a limiting process $\omega$ on bounded continuous 
functions $f(t)$ of $t\in [0,\infty]$, where $\lim_\omega f$ neglects all functions which tend to $0$ at $\infty$. One lets  $\alpha>0$ and defines the notation
\begin{equation}\label{power}
\lim_{\omega^\alpha} f:=\lim_\omega g, \ \ \  g(u):=f(u^\alpha)\qqq u\geq 0
\end{equation}
Then the relevant additional requirement is power invariance, \ie the equality of $\lim_{\omega^\alpha}$ with $\lim_{\omega}$. It can always be achieved due to the amenability of solvable groups and is technically very useful as shown in \S 7 of \cite{CSZ}. We refer to \cite{LSZ1} for a survey of the many key  developments around the singular traces whose prototype examples were uncovered by Dixmier.

\subsubsection{Wodzicki residue}

The Wodzicki residue \cite{MW} is a remarkable discovery which considerably enhances the power of the quantized calculus because it allows one, under suitable assumptions, to give a sound meaning to the integral of infinitesimals whose order is strictly less than one, while it agrees (up to normalization) with the Dixmier trace for infinitesimals of order $\geq 1$. Just as the Dixmier trace, it vanishes for infinitesimals of order $>1$ but  the extension of its domain is a considerable improvement which has no classical meaning. A striking example is given by the following new line in the dictionary (up to normalization):
\vspace{0.5cm}
 \begin{mdframed}[backgroundcolor=yellow!3]
\begin{tabular}{r | l}
& \\
 Scalar curvature of manifold of dimension $4$ & Functional $\cA\ni f\mapsto \displaystyle{\int\!\!\!\!\!\! -} f ds^2=\displaystyle{\int\!\!\!\!\!\! -}R f ds^4$ \\
&
\end{tabular}
\end{mdframed}

\vspace{0.5cm}
The Wodzicki residue  was extended to the framework of spectral triples with simple dimension spectrum in our joint work \cite{cmindex} with Henri Moscovici where we obtained a local expression (local in the above sense) for the cyclic cocycle which computes the index pairing (see \S \ref{sectcyclic}). There is a revealing example \cite{acsu2} where the computation of this local formula displays well the simplifications which occur when on manipulates infinitesimals (in the operator theoretic sense) modulo those of high order. The explicit form of the operators provided by the quantum group is quite involved but when one works modulo higher order, all sorts of irrelevant details are wiped out. 
The noncommutative residue was extended by R. Ponge to Heisenberg manifolds in \cite{Ponge1}
who obtained new invariants of CR-geometry  in a series of papers \cite{Ponge2}.

\subsubsection{Locality and  spectral action} \label{sectlocalact}
Another striking result obtained since the year 2000, was obtained in our joint work with Ali Chamseddine \cite{cc4}. It gives a general local formula for the scale independent terms $\zeta_{D+A}(0)-\zeta_{D}(0)$ in the spectral action, under the perturbation of $D$ given by the addition of a gauge potential $A$. Here  
$(\cA,\cH,D)$ is a spectral triple with simple dimension spectrum consisting of positive
integers and the general formula has finitely many terms and gives
\begin{equation}\label{zetaDplusA}
\zeta_{D+A}(0)-\zeta_{D}(0)=\,-\,{\int\!\!\!\!\!\! -}\, \log
(1+\,A\,D^{-1})=\,\sum_n\,\frac{(-1)^{n}}{n}\,{\int\!\!\!\!\!\! -}\,
(A\,D^{-1})^{n} .
\end{equation}
In dimension $4$ this  formula reduces to 
\begin{equation}\label{counfergrzeta} \begin{array}{rl}
\zeta_{D+A}(0)-\zeta_{D}(0)= & - \displaystyle{\cutint} AD^{-1}+\frac{1}{2}
 \displaystyle{\cutint}   (AD^{-1})^{2} \\[3mm] & -\frac{1}{3} \displaystyle{\cutint}
 (AD^{-1})^{3}+\frac{1}{4}\displaystyle{\cutint} (AD^{-1})^{4} . \end{array}
\end{equation}
and under the hypothesis of vanishing of the tadpole
term,  the variation \eqref{counfergrzeta}
is the sum of a Yang--Mills action and a
Chern--Simons 
 action relative to a cyclic $3$-cocycle
on the algebra $\cA$. This result should be seen as a first step in the long term program of performing the renormalization technique entirely in the framework of spectral triples. We refer to  \cite{walter1,walter2,walter3} for fundamental steps in this direction. After the initial breakthrough discovery in our joint work with D. Kreimer of the conceptual meaning of perturbative renormalization as the Birkhoff decomposition, we unveiled in our joint work with M. Marcolli \cite{CMbook} Theorem 1.106 and Corollary 1.107, an affine group scheme which acts on the coupling constants of physical theories as the fundamental ambiguity inherent to the renormalization process. 

At the conceptual level one can expect that since the local functionals filter out the quantum details, and give a semiclassical picture of the quantum reality, the knowledge of this semiclassical picture only allows to guess the quantum reality. In fact the interpretation of the line element $ds$ as the fermion propagator (in the Euclidean signature) already gives a hint of the modifications of geometry due to the quantum corrections. Indeed these corrections ``dress" the Fermion propagator as shown in Figure \ref{dsdr} :
\begin{figure}[H]
\begin{center}
\includegraphics[scale=0.5]{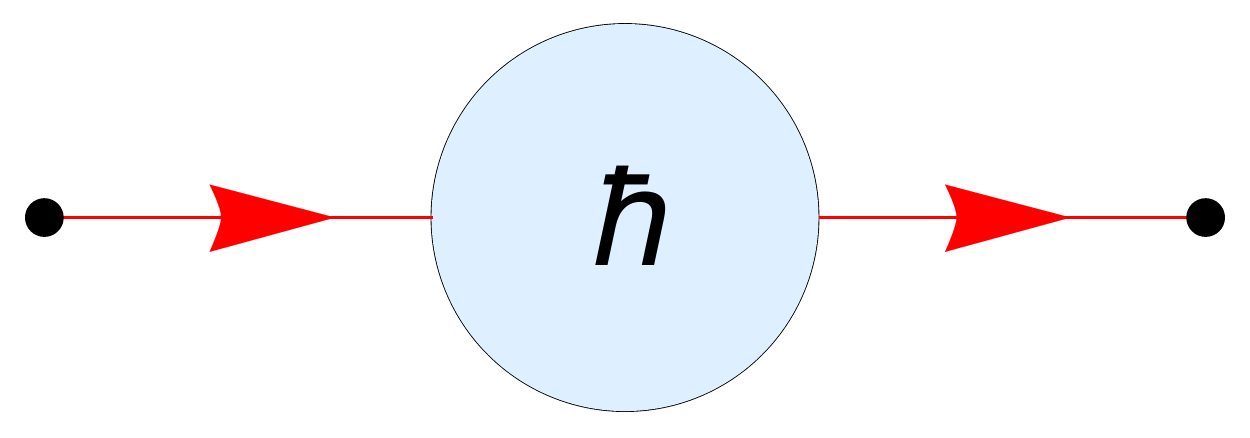}
\end{center}
\caption{Dressed Propagator \label{dsdr} }
\end{figure}
This dressing of the propagator should be interpreted physically as quantum corrections to the geometry, \ie the replacement of the initial ``classical" spectral triple by a new one in which the inverse line element undergoes the same dressing. In perturbative renormalization, this can only be expressed in the form of an asymptotic series in powers of $\hbar$ which modifies $D$ by multiplication by such a series of functions of the operator $D^2/\Lambda^2$. This suggests that in reality the geometry of space-time is far more interrelated to the quantum and should not be fixed and frozen to its semi-classical limit when one performs renormalization. The geometric paradigm of noncommutative geometry is well suited for this interaction with the quantum world. We refer to the books \cite{walterbook} of W. van Suijlekom and \cite{Iochumbook} of M. Eckstein and B. Iochum for complete accounts of the above topics on the spectral action.

\section{Cyclic cohomology}\label{sectcyclic}

The quantized calculus briefly presented in the dictionary replaces the differential $df$ of a function by the commutator $[F,f]$ which is an operator in Hilbert space where the self-adjoint operator $F$ acts and fulfills $F^2=1$. One can then develop the calculus using products of $k$ one forms as forms of degree $k$ and the graded commutator with $F$ as the differential. One has $d^2=0$ because $F^2=1$. Moreover under good summability conditions the higher forms will be trace class operators and puting things together one will get an abstract cycle of degree $n$ on the algebra $\cA$ \ie a graded differential algebra $\Omega$  which has $\cA$ in degree zero and is equipped with a closed graded trace $\tau: \Omega^n\to \C$. The integer $n$ is only relevant by its parity and by being large enough so that all $n$ forms are trace class operators. In 1981 (see 
the  Oberwolfach 1981 report  \cite{[Co$_{15}$]} on ``Spectral sequence and homology of currents for operator algebras'')  I discovered  cyclic cohomology,   the SBI
sequence and the related spectral sequence from the above framework. In particular the periodicity operator $S$ in cyclic cohomology was imposed by the ability to replace $n$ by $n+2$ in this framework. Since then cyclic cohomology plays a central role in noncommutative geometry as the replacement of de Rham cohomology. Moreover and parallel to what happens in differential geometry 
the integral pairing between $K$-homology and $K$-theory is computed by the pairing of the Chern characters in cyclic theory according to the diagram:

\vspace{0.5cm}

 {\bf\color{blue}
\begin{equation}\label{KHdiagram}
\xymatrix@C=45pt@R=45pt{
 K-\text{Theory} \ \ar[d]_{Ch_*}\ar@{<->}[rr] && K-\text{Homology} \ar@{->>}[d]^{Ch^*}
\\
HC_*\ \ar@{<->}[rr]&& HC^*}
\end{equation}}

\vspace{0.5cm}

The Chern character in $K$-homology is defined as explained above from the quantized calculus associated to a Fredholm module. However the obtained formula is non-local and difficult to compute explicitly. In \cite{Co-book} Theorem 8 of Chapter IV, I stated an explicit formula using the Dixmier trace for the Hochschild class of the Chern character, \ie its image under the forgetful map $I$ from cyclic to Hochschild cohomology. Following my unpublished notes a proof was included in the book \cite{FGV}. The result was then refined and improved in \cite{carey}.  

The Hochschild class of the character only gives partial information on the Chern character and 
a complete local formula, involving the Wodzicki residue in place of the Dixmier trace was obtained in 1995 in our joint work with H. Moscovici \cite{cmindex}. In fact index theory is an essential application of the cyclic theory and we shall briefly describe its role  in \S \ref{sectindex} below.
 
\vspace{0.5cm}
 \begin{mdframed}[backgroundcolor=yellow!3]
\begin{tabular}{r | l}
& \\
De-Rham homology of manifold  & Periodic Cyclic cohomology \\
&  \\
Atiyah-Singer Index Theorem & Local index formula in NCG \cite{cmindex}\\
&  \\
Characteristic classes and Lie algebra cohomology & Cyclic cohomology of Hopf algebras \\
& 
\end{tabular}
\end{mdframed}

\vspace{0.5cm}

\subsection{Transverse elliptic theory and cyclic cohomology of Hopf algebras}

 In our work with H. Moscovici on the transverse elliptic operators for foliations we were  led  to develop the theory of characteristic classes in the framework of noncommutative geometry as a replacement of its classical differential geometric ancestor. We also developed another key ingredient of the transverse geometry  which is the geometric analogue of the reduction from type III to type II and automorphisms and we used the theory of hypoelliptic operators to obtain a spectral triple describing the type II geometry.  In fact the explicit computations involved in the local index formula for this spectral triple dictated two essential new ingredients: first a Hopf algebra $\cH_n$ which governs  transverse geometry in codimension $n$, second a general construction of  the cyclic cohomology of Hopf algebras as a far reaching generalization of Lie algebra cohomology \cite{CMos,CMos2,CMos3} and \cite{AK,HKRS1,HKRS2}. This theory has been vigorously developed by H. Moscovici and B. Rangipour \cite{MR0,M1,MR1,MR2}, the DGA version of the Hopf cyclic cohomology and the characteristic map from Hopf cyclic to cyclic were investigated by A.~Gorokhovsky in \cite{Goro}. 
\subsection{Cyclic cohomology and archimedean cohomology}  
   
 In \cite{CC6} with C. Consani, we showed that cyclic homology which was invented for the needs of noncommutative geometry is  the right theory to obtain, using the $\lambda$-operations,  Serre's
archimedean local factors of the $L$-function of an arithmetic variety as regularized determinants. Cyclic homology provides a conceptual general construction of the sought for ``archimedean cohomology''. This resurgence of cyclic homology in the area of number theory was totally unexpected and is a witness of the coherence of the general line of thoughts. In fact the recent work \cite{CCMT,CCMT1} of  
G. Cortinas, J. Cuntz, R. Meyer and G. Tamme relates rigid cohomology to  cyclic homology.
\subsection{Topological cyclic homology}
We refer to the book \cite{DGM} for a complete view of the relation of topological Hochschild and cyclic theories with algebraic $K$-theory. At the conceptual level one may think of topological cyclic  homology as  cyclic homology performed over the smallest possible ground ring which in algebraic topology is not the familiar ring $\Z$ of integers but the far more absolute sphere spectrum. The great advantage of working over this base is that the analogue of the tensor product becomes much more manageable as far as the action of permutations on the factors as well as their fixed points are concerned. This fact gives rise to new operations which do not exist in the general framework. These operations include an analogue of the Frobenius operator and in his work with I. Madsen \cite{H2,HM} and then in \cite{H0,H1,H3}, L. Hesselholt  has shown that the local and global Witt construction, the de Rham-Witt complex and the Fontaine theory all emerge naturally  in this framework.   All the more since he showed  recently how topological periodic cyclic homology with its inverse Frobenius operator may be used to give a cohomological interpretation of the Hasse-Weil zeta function of a scheme smooth and proper over a finite field in the form entirely similar to \cite{CC6}. Finally, the paper \cite{CCsalg}   shows that topological cyclic homology is  cyclic homology when applied to the Gamma-rings of G. Segal and   that this framework, based on the sphere spectrum, subsumes all previous attempts of developing  ``absolute algebra". We refer to the work of B.~Dundas in \cite{Dundas} for a survey of the topological  cyclic theory and an attempt towards extending the construction in order to incorporate the unstable information.

\section{Overall picture, interaction with other fields}

\begin{figure}
\begin{center}
\includegraphics[scale=0.35]{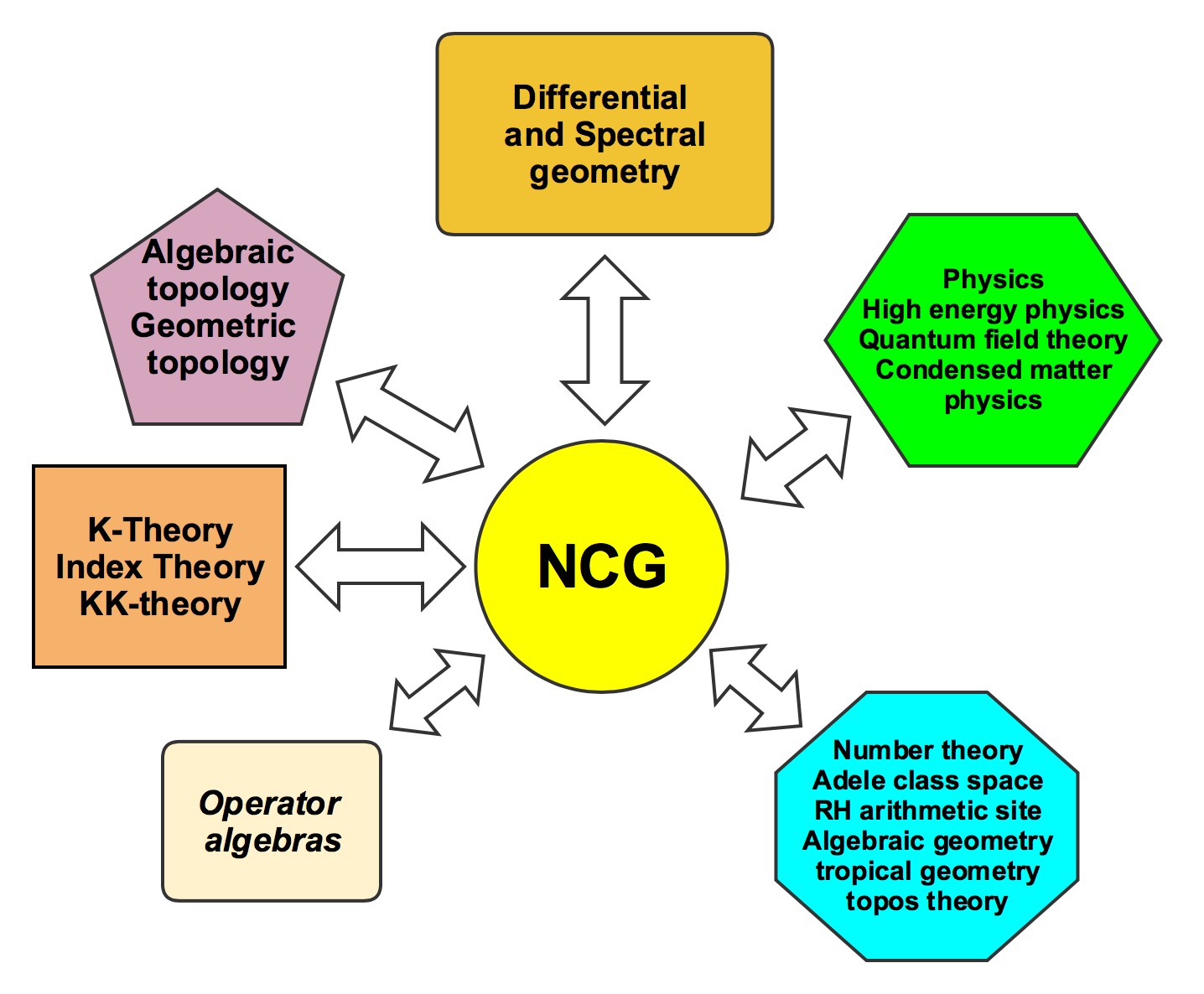}
\end{center}
\caption{Interactions with noncommutative geometry \label{mindmap1} }
\end{figure}

We describe in this section a few of the interactions of noncommutative geometry with other fields with no claim of being exhaustive.

\subsection{Index theory}\label{sectindex}

\subsubsection{Type II index theory}:\label{type2index}

Atiyah's index theorem for covering spaces \cite{Atiyah1} gave a striking application of the type II index theory which he used in order to define the index of an elliptic operator acting on the non-compact Galois covering of a smooth compact manifold. In Atiyah's theorem the  type II index  takes the same value as on the compact quotient and is hence integer valued. In \cite{Cofoliation} the type II index theorem for foliations displayed the role of the Ruelle-Sullivan current associated to a transverse measure and the index as well as the dimensions of the solution spaces are real valued and give striking examples of real dimensions as discovered by Murray and von-Neumann. These real dimensions are best understood at the intuitive level in the same way as densities of infinite sets along the leaves relate to cardinalities of finite sets and this point of view is extensively developed in \cite{Cofoliation}. 

\subsubsection{Higher index theory}\label{secthigherindex}

The need for a higher index theory is already visible starting with the longitudinal index theorem for foliations which made use of a transverse measure as explained in Section \ref{type2index}. Such transverse measures do not always exist and the theory would seem to be limited to the type II situation as opposed to the type III which often occurs in geometry.  Fortunately, the need of a transverse measure for the formulation of the theorem can be overcome when the result is formulated in terms of $K$-theory as was done in our joint work with G. Skandalis \cite{CS}. But in order to obtain numeric indices one still needs to define natural numerical invariants of $K$-theory. This is obtained using the pairing of $K$-theory with cyclic cohomology   and the early applications of the theory go back to the work on the transverse fundamental class \cite{transclass}. The strategy  can be divided in several steps
\begin{enumerate}
\item	Prove analytic invariance or vanishing properties of the analytic index as an element of the $K$-theory of a suitable algebra (most often a $C^*$-algebra associated to the geometric situation). 
\item Compute the pairing of the analytic index with cyclic cocycles by explicit geometric formulae, in the spirit of the original result of Atiyah and Singer.
\item Show that the algebraic pairing with the cyclic cocycle in fact makes sense at the level of the $K$-theory of the $C^*$-algebra.
\end{enumerate}
Note that steps 1 and 3 are of hard analytic nature while the step 2 is of a more algebraic and computational nature. The third step is usually the most difficult. 
The early applications of the higher index theory   included the vanishing of $\hat A$-genus for a foliated compact manifolds with 
leafwise positive scalar curvature \cite{transclass} and \cite{ZZ,BH}, 
the Novikov conjecture for Gelfand-Fuchs classes \cite{transclass}
and the  Novikov conjecture for hyperbolic groups (\!\cite{[Co-M$_1$]}   see \S \ref{sectnovikov}).  The proof of the Novikov conjecture for Gelfand-Fuchs classes illustrates the power of cyclic cohomology theory and
there is no other approach available to prove this result which inspired the recent proof by S.~Gong, J.~Wu and G.~Yu of the Novikov conjecture for discrete subgroups of the group of volume preserving diffeomorphisms  in \cite{GWY}. M.~Puschnigg has made key contributions \cite{Puschnigg,Puschnigg1,Puschnigg2} to step 3 and was able to use it to  prove the Kadison-Kaplansky conjecture (which states that the reduced group $C^*$-algebra does not contain non-trivial idempotents) for hyperbolic groups. Moreover R.~Meyer fully developed in \cite{Meyer-cyclic} the cyclic theory in the context of bornological algebras. A basic new idea is the notion of analytic nilpotence for bornological algebras which allows him to use the Cuntz-Quillen approach to cyclic theories. In \cite{MP,MP1}, H. Moriyoshi and P. Piazza proved an higher index analogue of the  Atiyah-Patodi-Singer formula for foliated manifolds using the pairing with the Godbillon-Vey cocycle. In \cite{HR3,HR4,HR5} N. Higson and J. Roe construct an analytical analogue of surgery theory and connect this theory with the usual surgery theory for manifolds. 
They define an analytic analogue of algebraic Poincar\' e complexes and show that such a complex has a signature, lying in the $K$-theory of the appropriate $C^*$-algebra, and that this signature is invariant under
bordism and homotopy. 
 They construct  a natural commutative diagram of exact sequences, sending the surgery exact sequence of Wall  to an analytic surgery sequence improving on previous results of J. Rosenberg \cite{Rosenberg}. 
The above gives only a brief glance at  the primary higher index theory and
we refer to the survey \cite{GE} of  E.~van Erp and A.~Gorokhovsky for an excellent  introduction to the subject  as well as a detailed account of the local index formula in NCG. 

We now turn to the recent developments \cite{PiS,XY0,XY1,XY2,MP,GMP,Goro} in the secondary higher index theory.
When the primary index vanishes as in the case of positive scalar curvature, or in the relative situation of homotopic manifolds for the signature, the secondary theory comes into play, and one obtains a secondary invariant called the rho-invariant.  This theory was first introduced by Higson and Roe \cite{HiR,HR5}. We refer to the splendid survey by Zhizhang Xie and Guoliang Yu  \cite{XY2} for the remarkable recent developments in this domain. 
 We shall not attempt (for lack of space) to do justice to these considerable  progress but briefly mention the papers  \cite{PiS} of P. Piazza and T. Schick and \cite{XY0,XY1,XY2} of  G. Yu and Z. Xie  who develop sophisticated secondary invariants involving the group $C^*$-algebra of Galois covers to analyse the connected components of    
   the space of metrics of positive scalar curvature.  They  
 proved in full generality  the connection of the higher index of the Dirac operator on a spin manifold with boundary to the higher rho invariant of the Dirac operator on the boundary, where the boundary is endowed with a positive scalar curvature metric. In \cite{WXY} with S. Weinberger, G. Yu and Z. Xie  settled a long standing question by showing that   the higher rho invariant is a group homomorphism. They  then apply this result to study non-rigidity of topological manifolds by giving a lower bound for the size of the reduced structure group of a closed oriented topological manifold, by the number of torsion elements in the fundamental group of the manifold.  Note that step 3 above comes into play when one uses cyclic cohomology to obtain numerical invariants by pairing with the rho invariant. In \cite{CWXY} the pairing of delocalized cyclic cocycles with the higher rho invariants is computed in terms of the delocalized eta invariants which were shown in  \cite{XY3} to be algebraic if the Baum-Connes conjecture holds, thus giving a new test of the conjecture. We refer to 
\cite{XY2} for an account of this theory.

\subsubsection{Algebraic index theory}

The algebraic index theory is a vast subject which has been successfully developed  by R. Nest and B. Tsygan (see \cite{NT1,NT2,NT3,NT4,NT5}) and with P.~Bressler and  A.~Gorokhovsky in \cite{BGNT}. It is an approach to the Atiyah-Singer index formula using deformation quantizations. It gives a powerful way of understanding conceptually and proving the algebraic index formulas which play for instance  a key role in step 2 of the higher index theory (\S \ref{secthigherindex}). The algebra of pseudo-differential operators provides a deformation quantization of the algebra of smooth functions on the cotangent space of a manifold. The index map is formulated using the trace density, which gives a map from the Hochschild or negative cyclic homology of the deformation quantization to the de Rham cohomology of the manifold. We refer to \cite{NT1,BGNT} for this  line of development as well as \cite{PPT,PPT1} for index theory in the framework of Lie groupoids and to \cite{Nistor,Monthubert,Monthubert1,MN} for the analysis on singular spaces  consisting of a smooth interior part and a tower of foliated structures at the boundary, including corner and cone points and edges of all possible dimensions.  This leads us to the next topic which is the geometrization of the pseudodifferential calculus and of deformations.

\subsection{Smooth groupoids and their $C^*$-algebras	in  geometry}\label{sectgroupoid}

In \cite{Co-book}  the construction of the tangent groupoid of a manifold was presented as a geometric framework for the theory of pseudo-differential operators.  The tangent groupoid embodies in a geometric manner the deformation from symbols to operators. The main tool there is the functorial association of a $C^*$-algebra $C^*(G)$ to a smooth groupoid $G$. As in the case of discrete groups, it admits a reduced and unreduced version, a distinction which disappears in the amenable case. This geometrization of operators allows one to give a completely geometric proof of the Atiyah-Singer index theorem as sketched in \cite{Co-book} and shown in \cite{Higson}. In a remarkable series of papers C. Debord and G. Skandalis \cite{DS1,DS2,DS3} have considerably extended this idea and advanced in a general program of geometrization of  the existing pseudo-differential calculi. In her work C. Debord had shown how to construct the holonomy groupoid for a singular foliation defined by a Lie algebroid  when the fiber dimension is the same as the leaf dimension. The general construction of the holonomy groupoid for a singular foliation was done by Androulidakis and Skandalis and this was the starting point of the index theory in this framework. In \cite{DS4}, C.Debord and G. Skandalis construct Lie groupoids by using the classical procedures of blowups and of deformations to the normal cone. They show that the blowup of a manifold sitting in a transverse way in the space of objects of a Lie groupoid leads to a calculus  similar to the Boutet de Monvel calculus for manifolds with boundary.  Their geometric constructions gives rise  to  extensions of $C^*$-algebras and they compute the corresponding K-theory maps and the associated index. Smooth groupoids and the associated $C^*$-algebras thus  provide a powerful tool in differential geometry. We refer to  \cite{DS5} for an excellent survey of this theory.  This tool also applies successfully in representation theory as shown by N. Higson in \cite{Higson1} who gave substance to the idea of G. Mackey on the analogy between complex semisimple groups and their Cartan motion groups.


\subsection{Novikov conjecture for hyperbolic groups}\label{sectnovikov}
In the C. R. Acad. Sc. Paris Note of 1988 ``Conjecture de Novikov et groupes hyperboliques''  \cite{[Co-M$_5$]} and in  \cite{[Co-M$_1$]},   Henri Moscovici and myself proved the Novikov conjecture for hyperbolic groups using cyclic 
cohomology combined with analytic techniques. The basic observation which was part of the origin of cyclic cohomology, is that, given a discrete group $\Gamma$, a group cocycle $c\in Z^n(\Gamma,\C)$ gives a cyclic cocycle $\tau_c$ on the noncommutative group ring $\C[\Gamma]$. Moreover an index formula \cite{[Co-M$_5$]}  \cite{[Co-M$_1$]} shows that the pairing of the cyclic cocycle gives the higher signature. The main difficulty is to pass from the very small domain of definition of the cocycle provided by $\C[\Gamma]$ to a larger domain having the same $K$-theory as the $C^*$-algebra $C_r^*(\Gamma)$ of the group, in which the signature is known to be homotopy invariant. The key technical tool which we used is the Haagerup property which provides a good notion of Schwartz rapid decay elements in $C_r^*(\Gamma)$ for hyperbolic groups. This together with the existence of bounded group cocycles representing cohomology classes of dimension $>1$ (which is dual to the non-vanishing of the Gromov norm) gave the solution of step 3 of the higher index technique of \S \ref{secthigherindex} and proved  the Novikov conjecture for hyperbolic groups \cite{[Co-M$_1$]}. Subsequently, and together with M. Gromov and H. Moscovici we produced another proof based on geometry. The original proof based on analysis and subalgebras of $C^*$ algebras stable under holomorphic functional calculus is still a main tool for   
 the Baum-Connes conjecture as in the work of V. Lafforgue who proved the latter conjecture for all hyperbolic groups (see \S \ref{sectbcconj}).

\subsection{The Baum-Connes conjecture} \label{sectbcconj}
In this section we give a short overview of the Baum-Connes conjecture and refer to \cite{GJV} for a remarkable survey of this topic. At the beginning of 1980, I had shown that closed transversals of foliations provide idempotents in the $C^*$-algebra of the foliation. This geometric construction played a key role in the explicit construction of the finite projective modules for the noncommutative torus \cite{C} as described in \S \ref{sectnctorus}. It was however quite clear that this geometric construction could not describe all of the $K$-theory of $C^*(V,F)$ even in the simplest case of a fibration. I met Paul Baum at the Kingston conference of $1980$ and he explained his work with Ron Douglas describing the $K$-theory of a manifold $W$ using geometric cycles given by triples $(M,E,f)$ where $M$ is a compact manifold, $E$  a vector bundle over $M$ and $f:M\to W$ is a continuous map, while one makes suitable $K$-orientation hypothesis. This construction immediately gave the clue on how to extend the construction of $K$-theory classes from transversals to a much more general and flexible one; the principle at work there is a general principle of noncommutative geometry which is worth emphasizing because of its potential use in other circumstances, it can be stated as follows: 

\vspace{0.5cm}
 \begin{mdframed}[backgroundcolor=yellow!3]
 While it is in general difficult to map a noncommutative space $X$ (such as a leaf space) to an ordinary space, it is quite easy to construct maps the other way, and use such maps to test homological properties of $X$.
 \end{mdframed}
 
 \vspace{0.5cm}
 
 In fact, in the case of foliations $(V,F)$, in order to get a map $f:M\to V/F$ from a manifold to a leaf space it is enough to cover $M$ by open sets and maps $f_j:\Omega_j\to V$ and to provide a one-cocycle 
ensuring that the maps match in the leaf space. 
The Baum-Connes conjecture was formulated by P. Baum and myself in 1981 and had a great impact on the K-theory of operator algebras (see \cite{GJV}). It was first stated for discrete groups and foliations, for Lie groups (as the Connes-Kasparov conjecture) and then  generalized (see Baum--Connes--Higson \cite{BCH}, and J-L. Tu \cite{Tu}) to the framework of crossed products by locally compact groupoids. It has been a major application of the Kasparov bivariant theory. By its relevance for representation theory of real or $p$-adic Lie groups, the analysis of  discrete groups, quantum groups, group actions etc,  it is one of the most powerful and unifying principles emerging as a link between algebraic topology and analysis.   It holds in full generality for Lie groups and in the amenable framework, and very powerful results have been obtained (see e.g.  \cite{HGT}) but its limits were found in \cite{HLS}  where counterexamples were obtained.

One of the most remarkable achievements since 2000 is the breakthrough contributions of Vincent Lafforgue who developed the analogue of the Kasparov bivariant theory in the context of Banach algebras \cite {Vlafforgue0} and successfully applied this new tool combined with his work on the comparison of $K$-theory in the Banach and $C^*$-contexts. He was able to cross the obstruction given, for discrete groups, by the property T of Kazhdan which had from the beginning in the 1980's forced to use the reduced $C^*$-algebra of the group as the natural one for the general formulation, but seemed an insurmountable obstacle in general. In \cite{Vlafforgue_these} V. Lafforgue establishes the Baum-Connes conjecture for many groups with property $T$. In \cite{MY}  I.~Mineyev and G.~Yu showed that hyperbolic groups are strongly bolic in the sense of V. Lafforgue and therefore satisfy the Baum-Connes conjecture. The results of Lafforgue cover not only reductive real and p--adic Lie groups but also  discrete cocompact subgroups of real Lie groups of real rank one and  discrete cocompact subgroups of real Lie groups  $SL_3$ over a local field. The conjecture with coefficients was proved by P. Julg for the property $T$ group $Sp(n,1)$ in \cite{Julg}. Finally in  \cite{Vlafforgue} V. Lafforgue proved the Baum-Connes conjecture with coefficients for hyperbolic groups, a stunning positive result. 

\subsection{Coarse geometry and the coarse BC-conjecture} \label{sectcoarse}

Coarse geometry is the study of spaces from a ``large-scale" point of view. Its link with $C^*$-algebras and index theory is due to John Roe \cite{Roe,Roe1,Roe2,Roe4}. The book \cite{Roe} gives a self-contained and excellent introduction to the topic. 
In his thesis John Roe managed to develop the index theory for non-compact spaces such as the leaves of a foliation of a compact manifold without the need of the ambient compact manifold (which was needed for the $L^2$ index theorem for foliations) and he associated a $C^*$-algebra to a coarse geometry thus allowing to apply the operator technique in this new context. His construction is obtained by abstracting the properties which in the foliation context hold for  operators along the leaves coming from the holonomy groupoid. Let $(X,d)$ be a  metric space in which every closed ball is compact. Let $\cH$ be a separable Hilbert space equipped with a faithful and non-degenerate representation of $C_0(X)$ whose range contains no nonzero compact operator. The Roe algebra associated to a coarse geometry is the $C^*$-algebra $C^*(X,d)$ which is the norm closure of  operators $T$ in $\cH$ which are  locally compact\footnote{\ie such that $fT$ and $Tf$ are compact operators for any $f\in C_0(X)$} and of finite propagation\footnote{\ie their support stays at a finite distance from the diagonal. A pair $(x,y)\in X\times X$ is not  in the support of $T$ iff there are $f,g\in C_0(X)$, $f(x)\neq 0$, $g(y)\neq 0$, $fTg=0$. }, two notions which have a clear geometric meaning using the metric $d$. This $C^*$-algebra is independent of the choice of the representation of $C_0(X)$.  The BC-conjecture is adapted to the coarse context under the hypothesis that $(X,d)$ has bounded geometry and using a discrete model.

In \cite{Yu1} Guoliang Yu proved the Novikov conjecture for groups with finite asymptotic dimension, and in \cite{Yu2} he proved the coarse BC-conjecture for metric spaces which admit a  uniform embedding into Hilbert space, and as a corollary the Novikov conjecture for discrete groups whose coarse geometry is embeddable in Hilbert space. Moreover he introduced  Property A for metric spaces and showed that it implies uniform embedding into Hilbert space. We refer to the textbook  \cite{NY} for large scale geometry and the many other applications of property A.  

At the conceptual level, the coarse point of view is best expressed in Gromov's words 

\begin{small}
\begin{quote}
To regain the geometric perspective one has to change one's position and move the observation point far away from X. Then the metric in X seen from the distance $\delta$ becomes the original distance divided by $\delta$ and as $\delta$ tends to infinity the points in X coalesce into a connected continuous solid unity which occupies the visual horizon without any gaps or holes and fills our geometers heart with joy....
\end{quote}
\end{small}

This point of view is clearly fitting with the principle of locality in noncommutative geometry explained in \S \ref{sectlocal} provided one works ``in Fourier" \ie one deals with momentum space. The large scale geometry of momentum space corresponds to the local structure of space. Thus the momentum space of a spectral geometry $(\cA,\cH,D)$  should be viewed as a coarse space in a suitable manner.  A first attempt is to consider the distance $d(A,A'):=\Vert A-A'\Vert$ defined by the operator norm on the gauge potentials, \ie the self-adjoint elements of the $\cA$-bimodule of formal one forms $\Omega^1:=\{\sum a_jdb_j\mid a_j,b_j\in \cA\}$ which one endows with 
\begin{equation}\label{forms}
d(A,A')=\Vert A-A'\Vert_D, \ \ \text{where}\ \ \Vert a_jdb_j \Vert_D:= \Vert \sum a_j[D,b_j]\Vert
\end{equation}
The $\cA$-bimodule  structure should be combined with the metric structure to investigate this coarse space. Note that the ``zooming out" for the coarse metric corresponds to the rescaling $D\mapsto D/\Lambda$ for large  $\Lambda$, which governs the spectral action. Under this rescaling $D\mapsto D/\Lambda$, 
the metric $d(\phi,\psi)$ on states of \eqref{dirac distance1} is multiplied by $\Lambda$, \ie one zooms in, which corresponds to the duality between the two metric structures, the first given by \eqref{dirac distance1} on the state space and the second given by \eqref{forms} on gauge potentials.     
\subsection{Comparison with toposes}\label{sectcomptopos}

Both noncommutative geometry and the theory of toposes provide a solution to the problem of the coexistence of the continuous and the discrete. The first by the use of the quantum formalism for real variables as self-adjoint operators and the second by a far reaching generalization of the notion of space which embodies usual topological spaces on the same footing as combinatorical datas such as small categories. 
The relation between  noncommutative geometry and the theory of toposes  belongs to the same principle as the relation established by the Langlands correspondence between analysis of complex representations on one hand and Galois representations on the other.

The functorial construction of the $C^*$-algebra associated to a smooth groupoid (\S \ref{sectgroupoid}) is the prototype of the correspondence between  noncommutative geometry and the theory of toposes. The role of the operator algebra side is to obtain  global results such as index theory,  homotopy invariance of the signature or vanishing theorems of the index of Dirac operators in the presence of positive scalar curvature. Such results would be inaccessible in a purely local theory.

  A systematic relation between the two theories has been developed by S. Henry in \cite{SimonH} and is based on the use of groupoids on both sides. 
  
  But as explained in  \cite{CaraL} the theory of Grothendieck toposes goes very far beyond its geometric role as a generalization of the notion of  space. The new input comes from its  relation with logics and the key notion of the classifying topos for a geometric theory of first order. We refer to \cite{CaraL} for an introduction to this fundamental aspect of the theory and the role of the duality between a topos and its presentation from a specific site. 
  
    From my own point of view, the great surprise was that very simple toposes provide, through their sets of points, very natural  examples of noncommutative spaces. The origin of this finding came from realizing that the classifying space of a small category $\cC$ gets much improved when it is replaced by the topos $\hat \cC$ of contravariant functors from $\cC$ to the category of sets.  In the C.R. Acad. Sc. Paris Note of 1983 ``Cohomologie cyclique et foncteurs Ext$^n$'' (\cf  \cite{CoExt}) I introduced the  cyclic category $\Lambda$ and showed that cyclic cohomology is a derived ext-functor after embedding the non additive category of algebras in the abelian category of $\Lambda$-modules. I also showed that the classifying space (in Quillen's sense) of the small category $\Lambda$ is the same as for the compact group $U(1)$ thus exhibiting the strong relation of $\Lambda$ with the circle. But since the classifying space of the ordinal category (of totally ordered finite sets and non-decreasing maps) is contractible, the information captured by the classifying space is only partial. This loss of information is completely repaired by the use of the topos $\hat \cC$ associated to a small category $\cC$. For instance the topos  thus associated to the ordinal category classifies the abstract intervals (totally ordered sets with a minimal and a maximal element) and encodes Hochschild cohomology as a derived functor. The lambda operations in cyclic homology \cite{Loday} enrich the cyclic category $\Lambda$ to the epicyclic category and the latter fibers over the category $(\star,\nt)$ with a single object and morphisms given by the semigroup $\nt$ of non-zero positive integers under multiplication. The topos $\wnt$ associated as above to this category underlies the arithmetic site (whose structure also involves a structure sheaf) and it was a great surprise to realize \cite{CCarith} that the points of this topos form a noncommutative space intimately related to the adele class space of \S \ref{sectzeta}. 
  
  One is surely quite far from a complete understanding of the inter-relations between topos theory and noncommutative geometry, and as another  
  instance of an unexpected relation we mention that the above topos $\wnt$ underlying the arithmetic site qualifies for representing the point in noncommutative geometry. More precisely noncommutative geometry works with (separable) $C^*$-algebras up to Morita equivalence and in each class there is a unique (up to isomorphism) $A$ which is stable \ie such that $A\simeq A\otimes \cK$  where $\cK$ is the algebra  of compact operators. The key fact then (\!\cite{CCarithgeo}, \S 8) is that the algebra $\cK$  of compact operators is naturally a $C^*$-algebra in the above topos, \ie it admits a natural action (unique up to inner) of the semi-group $\nt$ by endomorphisms which is hence inherited by any stable $C^*$-algebra.  Moreover the classification of matroids by J. Dixmier \cite{dix2} is given by the same space as the space of points of the topos and the corresponding algebras are obtained as stalks of the sheaf associated to $\cK$.

\subsection{Quantum field theory on noncommutative spaces} 
Noncommutative tori appear naturally in compactifications of string theory as shown in [140] (based on [141]), where the need for developing quantum field theory in such noncommutative spaces became apparent. The simplest case to consider is the Moyal deformation of Euclidean space. What became quickly discovered is that there is a mixing between ultraviolet and infrared divergencies so that the lack of compactness of the resolvent of the Laplacian becomes a real problem. In a remarkable series of papers H. Grosse and R. Wulkenhaar \cite{Grosse:2004yu} together with V. Rivasseau and his collaborators \cite{Rivasseau:2005bh} managed to overcome this difficulty and to renormalize the theory by adding a quadratic term which replaces the Laplacian by the harmonic oscillator.
They then developped Euclidean $\phi^4$-quantum field theory on four-dimensional Moyal space with harmonic propagation even further than its commutative analogue. For instance they showed that in contrast with the commutative case there is no Landau ghost in this theory: it is asymptotically safe \cite{Disertori:2006nq}, leading to the tantalizing
prospect of a full non-perturbative construction, including sectors with Feynman graphs of all genera. 

Grosse and Wulkenhaar went on to extensively study and fully solve the planar sector of the theory \cite{Grosse:2012uv} (certainly its most
interesting part since it  contains all the ultraviolet divergencies). This sector  also
identifies with an infinite volume  limit of the model in a certain sense. 
Very surprisingly both Euclidean invariance and translation invariance are restored in this limit.
Building on the identities introduced in \cite{Disertori:2006nq} they were able to fully
solve all correlation functions of a model in terms of a single sector of the two-point function 
They went even further and gave strong numerical evidence \cite{Grosse:2014lxa} (and a recent proof in the case of the two point function of the $\phi^3$ theory,  in collaboration with A. Sako \cite{Grosse:2016qmk}) 
that  Osterwalder-Schrader positivity also holds for such models. 
In particular, given the specific problems  (such as Gribov ambiguities) which have plagued 
and delayed the construction of ordinary non-Abelian gauge theories,
the planar sector of the Grosse-Wulkenhaar theory 
has now become our best candidate for a
rigorous reconstruction of a non-trivial Wightman
theory in four dimensions. 

The study of the Grosse-Wulkenhaar model has been also inspirational for 
promising generalizations of matrix models and non-commutative field theories, namely the tensor models 
and tensor field theories discovered and ceveloped by R. Gurau, V. Rivasseau and collaborators \cite{tensors}. We refer to \cite{Wsurvey} for an extensive survey of the remarkable developments of quantum field theory on noncommutative spaces.

\subsection{Noncommutative geometry and solid state physics}

Since the pioneering work of Jean Bellissard on the quantum Hall effect and the deep relation which he unveiled as the noncommutative nature of the Brillouin zone (see \cite{Co-book} Chapter IV, 6) the relation between noncommutative geometry and solid state physics has been vigorously pursued. We refer to the book \cite{Prod3} and papers \cite{Prod1, Prod2}
for the recent progress as well as to the papers of Terry Loring on the finite dimensional manifestations of the K-theory obstructions and of C.~Bourne, A.~Carey, A.~Rennie \cite{BCR} for  the bulk-edge correspondence.

\section{Spectral realization of zeros of zeta and the scaling site}\label{sectzeta}
We describe in this section the main steps from
quantum statistical mechanical models arising in number theory, the  adele class space as a noncommutative space, to the scaling site, an object of algebraic geometry underlying the spectral realization of the zeros of the Riemann zeta function.

\subsection{Quantum statistical mechanics and number theory}\label{sectbcsystem}

Discrete groups $\Gamma$ provide very non-trivial examples of factors of type II$_1$ and their left regular representation given by the left action of $\Gamma$ in the Hilbert space $\ell^2(\Gamma)$ of square integrable functions on $\Gamma$ generates a type II$_1$ factor $R(\Gamma)$ as long as the non-trivial conjugacy classes of elements of $\Gamma$ are infinite. It was shown by Atiyah that the type II index theory relative to $R(\Gamma)$ can be used very successfully for Galois coverings of compact manifolds. To obtain type III factors from discrete groups, one considers the relative situation of pairs $(\Gamma,\Gamma_0)$, where $\Gamma_0\subset \Gamma$ is a subgroup which is almost normal inasmuch as the left action of $\Gamma_0$ on the coset space $\Gamma/\Gamma_0$ only has finite orbits. The left action of $\Gamma$ in the Hilbert space $\ell^2(\Gamma/\Gamma_0)$ then generates a von Neumann algebra which is no longer of type II in general and whose modular theory depends on the integer valued function $L(\gamma):=$ cardinality of the orbit $\Gamma_0(/\gamma\Gamma_0)$ in the coset space $\Gamma/\Gamma_0$. In fact the commutant von Neumann algebra is generated by the Hecke algebra $\cH( \Gamma,\Gamma_0)$ of convolution of functions $f$ of double cosets which have finite support in $\Gamma_0\backslash\Gamma/\Gamma_0$. The modular automorphism of the canonical state is given by \cite{bcsystem} 
$$
\sigma_t(f)=\left(\frac{L(\gamma)}{R(\gamma)}\right)^{it} f(\gamma), \ \ R(\gamma):=L(\gamma^{-1}).
$$
Almost normal subgroups $\Gamma_0\subset \Gamma$ arise naturally by considering the inclusion of points over $\Z$ inside points over $\Q$ for algebraic groups as in the construction of Hecke algebras. The BC system arises from the  ``$ax+b$'' algebraic group $P$.
By construction $P_\Z^+ \subset P_\Q^+$ is an inclusion $\Gamma_0
\subset \Gamma$ of countable groups, where $P_\Z^+$ and $P_\Q^+$
denote the restrictions to $a>0$. This inclusion fulfills the
above commensurability condition and the associated Hecke algebra with its dynamics is the BC system. Its main interest is that it exhibits spontaneous symmetry breaking \cite{bcsystem} when it is cooled down using the thermodynamics of noncommutative spaces (see \cite{CMbook} for this general notion) and that its partition function is the Riemann zeta function. Considerable progress has been done since 2000 on extending the BC-system to number fields. For lack of space we shall not describe these developments but  refer to \cite{CMbook} and to \cite{LNT,Larsen2}. The paper \cite{LNT} elucidates the functoriality of the construction  in Theorem 4.4, where they show that the construction of Bost-Connes type systems extends to a functor which to an inclusion of number fields $K\subset L$ assigns a $C^*$-correspondence which is equivariant with respect to their suitably rescaled natural dynamics. In the paper \cite{Larsen2} M. Laca, N. Larsen and S. Neshveyev succeeded to extend the properties of fabulous states of the BC system to arbitrary number fields. We refer to the introduction of their paper for an historical account of the various steps that led to the solution after the initial steps of Ha and Paugam. They consider the Hecke pair consisting of the group $P_K^+$
 of affine transformations of a number field $K$ that preserve the orientation in every real embedding and the subgroup $P_O^+$ of transformations with algebraic integer coefficients. They then proceed and show that the Hecke algebra $\cH$ associated to the pair $(P_K^+,P_O^+)$  fulfills perfect analogues of the properties of the BC--system including phase transitions and Galois action. More precisely they obtain an arithmetic subalgebra (through the isomorphism  of $\cH$ to a corner in the larger system established in \cite{LNT}) on which ground states exhibit the ``fabulous states" property with respect to an action of the Galois group ${\rm Gal}(K_{\rm ab}:H_+(K) )$, where $H_+(K)$ is the narrow Hilbert class field, and 
 $K_{\rm ab}$ is the maximal abelian extension of $K$.

We refer to \cite{Cuntz} for a survey of the important other line of developments  coming from the exploration by Joachim Cuntz and his collaborators of the $C^*$-algebra associated to the action of the multiplicative semigroup of a Dedekind ring on its additive group. Representations of such actions give rise to particularly intriguing problems and the study of the corresponding $C^*$-algebras has motivated many of the new methods and general results obtained in this area.

Another very interesting broad generalization of the BC-system has been developed by M.~Marcolli and G.~Tabuada in \cite{MT}. Using  the Tannakian formalism, they categorify the algebraic data used in constructing the system such as roots of unity, algebraic numbers, and Weil numbers. They study the  partition function, low temperature Gibbs states, and Galois action on zero-temperature states for the associated quantum statistical system and they show that in the particular case of the Weil numbers the partition function and the low temperature Gibbs states can be described as series of polylogarithms.

\subsection{Anosov foliation}\label{sectguillemin}

The factor of type III$_1$ generated by the regular representation of the BC-system is the same (unique injective factor of type III$_1$, \cite{haagerup}) as the factor associated to the Anosov foliation of the sphere bundle of a Riemann surface of genus $>1$. This fact suggested an  analogy between the geometry of the continuous decomposition (as a crossed product: type II$_\infty\rtimes \R_+^*$) in both cases which led me in \cite{aczetacr}, based on the paper of Guillemin   \cite{guillemin},  to the adele class space as a natural candidate for the 
trace formula interpretation of the explicit formulas of Riemann-Weil. It is worth explaining
briefly the framework\footnote{The same paper \cite{guillemin} was used two years later by Deninger in \cite{Deninger} to motivate his search for an hypothetical cohomology theory using foliations.} of   \cite{guillemin}. 
Let $\Gamma\subset SL(2,\R)$ be a discrete cocompact subgroup, and $V=SL(2,\R)/\Gamma$. Let $F$ be the foliation of $V$ whose leaves are the orbits of the action on the left of the subgroup $P\subset SL(2,\R)$ of upper triangular matrices. At the Lie algebra level $P$ is generated by the elements
$$
E^+=\left(
\begin{array}{cc}
 0 & 1 \\
 0 & 0 \\
\end{array}
\right),  \ \ H=\frac 12 \left(
\begin{array}{cc}
 1 & 0 \\
0 & -1 \\
\end{array}
\right)
$$
One denotes by $\eta$ and $\xi$ the corresponding vector fields. The associated flows are the horocycle and geodesic flows. They correspond to the one-parameter subgroups of $P\subset SL(2,\R)$ given by
$$
n(a)=\left(
\begin{array}{cc}
 1 & a \\
 0 & 1 \\
\end{array}
\right), \ \ g(t)=\left(
\begin{array}{cc}
 e^{t/2} & 0 \\
 0 & e^{-\frac{t}{2}} \\
\end{array}
\right)
$$
The horocycle flow is normalized by the geodesic flow, more precisely one has $g(t)n(a)g(-t)=n(a e^t)$. The von Neumann algebra of the foliation, $M=W(V,F)$ is a factor of type III$_1$ and its continuous decomposition is visible at the geometric level: the associated factor of type II$_\infty$ is the von Neumann algebra $N=W(V,\eta)$ of the horocycle foliation. The one parameter group of automorphisms $\theta_\lambda$ scaling the trace is given by the action of the geodesic flow by automorphisms of $N=W(V,\eta)$ which has clear geometric meaning by the naturality of the construction of the von Neumann algebra of a foliation. The trace on $N=W(V,\eta)$ corresponds to the Ruelle Sullivan current which is the contraction of the $SL(2,\R)$-invariant volume form of $M$ by the vector field $\eta$. The rescaling of the trace by the automorphisms $\theta_\lambda$ follows from  the rescaling of $\eta$ by the geodesic flow. In \cite{guillemin} a heuristic proof of the Selberg trace formula is given using the action of the geodesic flow on the horocycle foliation. Our interpretation of the explicit formulas as a trace formula involves in a similar manner the action of $\R_+^*$ on the geometric space associated to the type II$_\infty$ factor of the continuous decomposition for the BC-system. This geometric space is the adele class space (for the precise formulation we refer to \cite{CMbook}).


\subsection{The adele class space}
In \cite{Co-zeta} we gave a spectral realization of the zeros of the Riemann zeta function and of L-functions based on the above
action of the Idele class group on the noncommutative space of Ad\`ele classes which is the type II space associated to the BC system by the reduction from type III to type II and automorphisms.  This result determined a geometric interpretation of the Riemann-Weil explicit formulas as a Lefschetz formula and also a reformulation of the Riemann Hypothesis in terms of the validity of a trace formula. The understanding of the  geometric side of the trace formula  is simpler when one works with the full ad\`ele class space \ie one does not further divide by $\hatz$. The local contribution from a place $v$ of a global field $K$ is obtained as the distributional trace of an integral $\int f(\lambda^{-1})T_{\lambda}d^*\lambda $ of operators of the form, with $\lambda\in K_v^*$, 
$$
(T_\lambda\xi)(x):=\xi(\lambda x)\qqq x\in K_v
$$
where $K_v$ is the local completion of $K$ at the place $v$. The distributional trace of $T_\lambda$ is the integral $\int k(x,x)dx$ on the diagonal of the Schwartz kernel which represents $T_\lambda$ \ie by  
$$
(T_\lambda\xi)(x)=\int k(x,y)\xi(y)dy\Rightarrow k(x,y)=\delta(y-\lambda x)
$$
This gives for the trace the formula $\int k(x,x)dx=\int \delta(x-\lambda x)dx=\vert 1-\lambda\vert^{-1}$ using the change of variables $u=(1-\lambda)x$ and the local definition of the module which implies $dx= \vert 1-\lambda\vert^{-1}du$. For the convolution operator $\int f(\lambda^{-1})T_{\lambda}d^*\lambda $ this delivers the local contribution to the Riemann-Weil explicit formula as $\int_{K^*_v}f(\lambda^{-1})\vert 1-\lambda\vert^{-1}d^*\lambda$. The above computation is formal but as shown in \cite{Co-zeta} the more precise calculation even delivers the subtle finite parts involved in the Riemann-Weil explicit formula. In \cite{Meyer}, R. Meyer showed how, by relaxing the Sobolev condition, one can effectively reprove the explicit formulas as a trace formula on the ad\`ele class space.

\subsection{The scaling site}

\subsubsection{Hasse-Weil form of the Riemann zeta function}

 After a few years of direct attack using only analysis I came to the conclusion that a much better understanding of the geometry underlying the ad\`ele class space was required in order to make a vigorous advance toward the solution of this problem. This work was undertaken in the ongoing collaboration with C. Consani. Among the major results obtained in these past years is a geometric framework 
in which one can transpose to the case of RH many of the ingredients of the Weil proof for function fields as reformulated by Mattuck, Tate and Grothendieck. 
 The starting point \cite{CC0,CC1} is the determination and interpretation as intersection
number using the geometry of the ad\`ele class space, of the real counting function $N(q)$ which gives the complete Riemann zeta function by a Hasse-Weil formula in the limit $q\to 1$, in the line of the limit geometry on $\F_q$ for $q=1$  proposed by J.Tits and C.Soul\'e.

\subsubsection{Characteristic one} The limit $q\to 1$ is taken analytically in the above development, and this begs for an algebraic understanding of the meaning of ``characteristic one". Characteristic  $p$ is defined by the congruence $p\simeq 0$. While the congruence $1\simeq 0$ is fruitless, its variant given by $1+1\simeq 1$ opens the door to a whole field which we discuss in \cite{Crh} as the ``world of characteristic one" whose historical origin goes quite far back. It has strong connections with the field of optimization \cite{Gaubert, Gaubert1}, with tropical geometry, with lattice theory and very importantly it appeared independently 
 with the  school of ``dequantization",  of V. P. Maslov and his collaborators \cite{Maslov,Lit}.  They  developed a satisfactory algebraic framework which encodes the semiclassical limit of quantum mechanics and called it idempotent analysis.  The starting observation is that one can encode the limit $\hbar\to 0$ by simply conjugating the addition of numbers by the power operation $x\mapsto x^\epsilon$ and passing to the limit when $\epsilon\to 0$.   The new addition of positive real numbers is 
$$
\lim_{\epsilon \to 0}\left(x^{\frac 1\epsilon}+y^{\frac 1\epsilon}\right)^\epsilon=\max \{x,y\}=x\vee y
$$
When endowed with this operation as addition and with the usual multiplication, the positive real numbers become a semifield $\rmax$. It is of characteristic $1$, \ie $1\vee 1=1$ and contains the smallest semifield of characteristic $1$, namely the Boolean semifield $\B=\{0,1\}$. Moreover, $\rmax$ admits non-trivial automorphisms and one has 
$$
{\rm Gal}_\B(\rmax):=\Aut_\B(\rmax)=\R_+^*, \ \ \fr_\lambda(x)=x^\lambda \qqq x\in \rmax, \ \lambda \in \R_+^*
$$
thus providing a first glimpse of an answer to Weil's query in \cite{Weilcdc} of an algebraic framework in which the connected component of the idele class group would appear as a Galois group. The most striking discovery of this school of Maslov, Kolokolstov and Litvinov \cite{Maslov,Lit} is that the Legendre transform which plays a fundamental role in all of physics and in particular in 
 thermodynamics in the nineteenth century, is simply the Fourier transform  in the framework of  idempotent analysis.

 \subsubsection{Tropical geometry and Riemann-Roch theorems}

The tropical semiring $\N_{\rm min}=\N\cup\{\infty\}$ with the operations $\min$ and $+$ was introduced by Imre Simon in \cite{Simon} to solve a decidability problem in rational language theory. His work is at the origin of the term ``tropical" which was coined by Marco Schutzenberger. Tropical geometry  is a vast subject, see \eg \cite{Gelfand,  Mik,Sturm}. 
We refer to \cite{virotagaki} for an excellent introduction starting from the sixteenth Hilbert problem.  In its simplest form (\!\cite{GK}) a tropical curve is given by  a graph with a usual line metric on its edges. The natural structure sheaf on the graph is the sheaf  of real valued functions which are continuous, convex, piecewise affine with integral slopes. 
The operations on such functions are given by  $(f\vee g)(x)=f(x)\vee g(x)$ for all $x\in \Gamma$ and  the product  is given by point-wise addition. One also adjoins the constant $-\infty$ which plays the role of the zero element in the semirings of sections. One proceeds as in the classical case with the construction of the sheaf  of semifields of quotients and finds the same type of functions as above but no longer convex. 
 Cartier divisors make sense and one finds that the order of a section $f$  at a point $x\in\Gamma$ is given by the sum of the (integer valued) outgoing slopes. The conceptual explanation of why the discontinuities of the derivative should be interpreted as zeros or poles is due to Viro, \cite{viro} who showed that it follows automatically if one understands that, as seen when dealing with  valuations the sum $x\vee x$ of two equal terms  should be viewed as ambiguous with all values in the interval $[-\infty, x]$ on equal footing. In their work Baker and Norine \cite{BN} proved in the discrete set-up of graphs an analogue of the Riemann-Roch theorem whose essence is 
that  the inequality ${\rm Deg}(D)\geq g$ (where $g$ is the genus of the graph) for a divisor  implies that the divisor is equivalent to an effective divisor. Once translated in the language of the chip firing game (\opcit\!\!),  this fact is equivalent  to the existence of a winning strategy if one assumes that the total sum of dollars attributed to the vertices of the graph is $\geq g$. We refer to \cite{GK} for the tropical curve version  of the Baker and Norine   theorem. 

\subsubsection{The topos, its points and structure sheaf}
The major discovery \cite{CCarith, CCarithgeo,CCscal} is that of the  scaling site, the topos of $\N^\times$-equivariant sheaves  on the Euclidean half-line, which is obtained by extension of scalars from the arithmetic site. The  points  of this topos form exactly the  sector of the ad\`ele class space involved in RH, 
and the action of $\R_+^*$ on the adele class space corresponds to the action of the Frobenius.  This provides in full the missing geometric structure on the ad\`ele class space which becomes a tropical curve since the topos inherits, from its construction by extension of scalars, a natural sheaf  of regular functions as piecewise affine convex functions. This structure is central in the well known results on the localization of zeros of analytic functions which involve Newton polygons in the non-archimedean case and the Jensen's formula in the complex case. The new feature given by the action of $\N^\times$ corresponds to the transformation $f(z)\mapsto f(z^n)$ on analytic functions. 
 The Newton polygons play a key role in the construction of the analogue of the Frobenius correspondences as a one parameter semigroup of correspondences already defined at the level of the arithmetic site. 
Finally the restriction   to the periodic orbit of the scaling flow associated to each prime $p$  gives a quasi-tropical structure which turns this orbit into a variant $C_p=\R_+^*/p^\Z$ of the classical Jacobi description $\C^*/q^\Z$ of an elliptic curve. On $C_p$, \cite{CCscal1} develops the theory of Cartier divisors, determines the structure of the quotient  of the abelian group of divisors by the subgroup of principal divisors, develops the theory of theta functions, and proves the Riemann-Roch formula which involves real valued dimensions, as in the type II index theory. The current situation concerning the evolution of this strategy towards RH is summarized in the essay \cite{Crh} in the volume on open problems in mathematics, initiated by John Nash and  recently published by Springer. After developing homological algebra in characteristic one in \cite{CChomol} we understood in \cite{CCsurvey} how to go back and forth to the complex situation using the above Jensen formula so that the fundamental structure  takes place over $\C$ and  fits well with noncommutative geometry.

\end{document}